\documentclass[letterpaper,journal]{IEEEtran}

\usepackage{cite}
\usepackage{amsmath,amssymb,amsfonts}
\usepackage{algorithm,algpseudocode}
\usepackage{graphicx}
\usepackage{bm}
\usepackage{textcomp}
\usepackage{stfloats}
\usepackage{url}
\usepackage[english]{babel}
\usepackage{xcolor}
\usepackage{bbm}
\usepackage{amsthm}
\usepackage{amsmath}
\usepackage{enumitem}
\usepackage{tikz}
\usepackage[normalem]{ulem}
\usepackage{hyperref}
\hypersetup{
    colorlinks=true,   	
    linkcolor=red,      
    citecolor = [rgb]{0 0.7 0},   	
    filecolor=magenta, 	
    urlcolor=blue
}
\usepackage{adjustbox}

\def\BibTeX{{\rm B\kern-.05em{\sc i\kern-.025em b}\kern-.08em
    T\kern-.1667em\lower.7ex\hbox{E}\kern-.125emX}}
    
\newtheorem{theorem}{Theorem}
\newtheorem{proposition}{Proposition}
\newtheorem{lemma}{Lemma}

\def\transpose{\intercal}

\definecolor{dgreen}{rgb}{0.2,0.7,0.2}
\def\clC{{\cal C}}
\def\clG{{\cal G}}
\def\RL{{\mathbb R}}

\def\BBE{\mathbb E}

\def\clP{\cal P}
\def\clQ{\cal Q}

\def\bchi{\boldsymbol{\chi}} 
\def\hbchi{\hat{\boldsymbol{\chi}}}

\newcommand{\Xn}{{\mathbf{X}^n}}

\newcommand{\ww}{\mathbf{w}}

\newcommand{\IND}{{\mathbb I}}
\newcommand{\argmin}{\mathop{\rm arg\, min}}

\title{Bayesian causal discovery:\\ Posterior concentration and optimal detection
\thanks{This work has been funded in part by one or more of the following grants: ARO W911NF1910269,
ARO W911NF2410094,
DOE DE-SC0021417,
Swedish Research Council 2018-04359,
NSF CCF-2008927,
NSF RINGS-2148313,
NSF CCF-2200221,
NSF CCF-2311653,
ONR 503400-78050,
ONR N00014-22-1-2363,
NSF A22-2666-S003, and the
NSF Center for Pandemic Insights DBI-2412522
}
}

\author{%
Valentinian Lungu,~\IEEEmembership{Graduate Student Member, IEEE,}
    \thanks{Statistical Laboratory, Centre for Mathematical Sciences, 
	University of Cambridge, Wilberforce Road, Cambridge CB3 0WB, 
	UK. Email: \texttt{vml26@cam.ac.uk}. Supported by the 
	Heilbronn Institute for Mathematical Research and G-Research.}
\and
Joni Shaska,~\IEEEmembership{Graduate Student Member, IEEE,}
    \thanks{University of Southern California, Los Angeles,
    CA 90089, USA. Email: \texttt{shaska@usc.edu}.
    }
\and Ioannis Kontoyiannis,~\IEEEmembership{Fellow, IEEE,}
	\thanks{Statistical Laboratory, Centre for 
	Mathematical Sciences, University of Cambridge, Wilberforce Road, 
	Cambridge CB3 0WB, UK. Email: \texttt{ik355@cam.ac.uk}.
    Supported in part
    by the EPSRC-funded INFORMED-AI project EP/Y028732/1.}
\and    
Urbashi Mitra,~\IEEEmembership{Fellow, IEEE}
    \thanks{University of Southern California, Los
    Angeles, CA 90089, USA. Email: \texttt{ubli@usc.edu}.}
}

\date{\today}

\begin{document}

\maketitle

\maketitle

\thispagestyle{plain}
\pagestyle{plain}

\begin{abstract}
We consider the problem of Bayesian causal discovery
for the standard model of linear structural equations
with equivariant Gaussian noise. A uniform prior is placed on
the  space of  directed acyclic graphs
(DAGs) over a fixed set of variables and, given the graph, independent
Gaussian priors are placed on the associated linear coefficients
of pairwise interactions. We show that the rate at which the posterior
on  model space concentrates on the true underlying
DAG depends critically on its nature: If it is {\em maximal}, in the
sense that adding any one new edge would violate acyclicity,
then its posterior probability converges to~1 exponentially
fast (almost surely) in the sample size~$n$.
Otherwise, it converges at 
a rate no faster than $1/\sqrt{n}$.
This sharp dichotomy is 
an instance of the important general phenomenon that avoiding overfitting
is significantly harder than identifying all of the structure
that is present in the model. 
We also draw a new 
connection between the posterior distribution 
on model space and recent results on optimal
hypothesis testing in the related problem of edge detection. 
Our theoretical findings
are illustrated 
empirically through simulation experiments.    
\end{abstract}



\section{Introduction}

Understanding the causes of phenomena influenced by multiple variables 
can often be aided by models that represent causal relations via 
directed acyclic graphs (DAGs)~\cite{pearl2009causality}. 
{\em Causal discovery}, or
{\em structure learning}, is the identification of the true underlying 
DAG or of other appropriate structures. 
Applications of causal discovery
span a broad range of different areas, 
including, among many others, the problem of topology inference in wireless 
networks~\cite{testi2020blind}, the analysis of gene networks~\cite{THOR}, 
the evaluation of medical
treatments~\cite{wu:24}.
Applications have motivated extensive research resulting in a variety 
of causal discovery tools, including constraint-based 
methods~\cite{PCalgorithm}, score-based 
approaches~\cite{Peters_2013,lasso,NOTEARS}, and Bayesian methods~\cite{bayesdag,BCDNets,beingBayesian,annadani2021variational,dynamicbayesian,cao:19}. 
In this work, we adopt a Bayesian approach,
in the context of the standard causal model based 
on linear structural equations.
This is partly motivated by the fact that the
Bayesian setting
naturally provides tools for quantifying confidence 
levels for discovered structures, 
typically
in the form of posterior probabilities. Bayesian approaches have also been successful in quantifying epistemic uncertainty, facilitating downstream tasks such as causal effect estimation~\cite{toth2022}.

The  various methods which have been developed for quantifying uncertainty in causal discovery
include a bootstrap approach~\cite{BootStrap}; 
confidence intervals~\cite{chickering2002optimal,StriederDrtonconfidence};
and the formulation of causal discovery as a series of composite Neyman-Pearson hypothesis tests~\cite{NPCausalISIT}.
 In addition, several effective inference techniques have been employed for the setting considered in this work, such as a Gibbs sampler~\cite{heckerman2006a};
an efficient Metropolis-Hastings algorithm~\cite{viinikka2020scalablebayesianlearningcausal};
and a variational inference method for sampling node permutation matrices~\cite{BCDNets}.

The main difficulty in the practical
application of Bayesian methods is the high computational 
complexity of standard techniques,
especially in high-dimensions.
In order
to improve mixing in Markov chain Monte Carlo (MCMC) 
samplers, 
an `edge reversal' move is introduced in~\cite{edgereversal};
a sampler on the space of causal orderings rather than 
causal networks is used in~\cite{beingBayesian};
and a combination of heuristics and pruning
is used in~\cite{beingBayesian}. 
Variational methods have also been extensively
used~\cite{annadani2021variational,bayesdag,BCDNets,lorch2021dibs}. 

This work focuses on the linear causal model 
with homoscedastic additive Gaussian noise. 
We place a uniform prior on the space
of all relevant DAGs and independent Gaussian priors
on the {\em edge weights}, that is, 
on the linear coefficients of pairwise interactions.
Both 
classical~\cite{StriederDrtonconfidence,Chen_2019_ordering}
and Bayesian~\cite{BCDNets,bayesian_challenges,chang:24}
versions of this model have been considered,
and it is known to be identifiable
under mild conditions~\cite{Peters_2013}. 

Our main contributions are the following:

{\bf (i)}
	Let $\pi(S|\Xn)$ denote the posterior
	probability of a  DAG $S$.
	We show that, when 
	the observations $\Xn=(X_1,\ldots,X_n)$
	are generated under $S^*$,
	the rate at which 
    $\pi(S^*|\mathbf{X}^n)\to 1$ heavily depends on 
	$S^*$. If $S^*$ is {\em maximal}, i.e., if adding any one 
	new edge to it would violate the acyclicity constraint,
	then, as $n\to\infty$, $\pi(S^*|\mathbf{X}^n)\to 1$ exponentially fast
	(almost surely).
    	Otherwise, the convergence is much slower, 
	of rate no faster than
	$1/\sqrt{n}$.

{\bf (ii)}
	In the special case where the nonzero edge weights
	are {\em a priori} assumed to all be equal to~1, 
	we show that
	$\pi(S^*|\mathbf{X}^n)\to 1$
	exponentially fast
	(almost surely), with exponent
	equal to~$1/2$.

{\bf (iii)}
We demonstrate that Bayesian methods are optimal 
for edge detection under the Neyman-Pearson causal discovery 
framework of~\cite{NPCausalISIT}.

{\bf (iv)} We illustrate and validate the above theoretical results on simulated
	data. 

The sharp dichotomy displayed by the posterior convergence
results in~{\bf (i)}
can be viewed as an instance of 
the following general phenomenon:
Avoiding overfitting is significantly harder than identifying all the
structure that is actually present in the model. Here,
this means that, when $S^*$ is maximal, then
it is impossible to overfit, because it is impossible 
to add any spurious edges to the true model,
so
the posterior converges quickly.
By contrast,
when $S^*$ is not maximal, then much stronger evidence,
i.e., many more observations, 
are required to rule out all extra
edges, 
leading to 
much slower convergence. Notice that overfitting is not possible when all nonzero weights are fixed \textit{a prior} (such as in~{\bf (ii)}), since it only occurs when 
it is difficult to discern the absence of an edge 
from its presence with an arbitrarily small but 
still nonzero weight.
Earlier instances of this dichotomy have been noted,
e.g., in~\cite{Schwarz1978,kass:95,castillo:15,shin:18}.
In fact, as we will see in the proofs of our main
results in Sections~\ref{sec:binarycase}
and~\ref{sec:generalcase}, the analysis is
similar, at least in spirit,
to  Schwarz's derivation~\cite{Schwarz1978}
of the Bayesian information criterion (BIC).

A similar phenomenon was recently observed in the 
context of a different class of graphical models, 
namely, tree models for time series.
Specifically, for the Bayesian context tree (BCT)
model class~\cite{BCT-JRSSB:22,K-theory:24},
it was shown~\cite{papag-K-BA:24}
that the posterior on the model space converges
at a polynomial rate, unless the true underlying
model is the unique maximal
model, in which case the posterior convergence 
rate is exponential.

\section{The setting}
\label{sec:intro}

\subsection{The Bayesian linear causal model}

We consider the linear structural equations model
\begin{equation}
X=AX+\epsilon.
\label{eq:LSE}
\end{equation}
Here, $X=(X(1),\ldots,X(d))^\transpose\in\RL^d$ 
is a $d$-dimensional sample, and
$\epsilon=(\epsilon(1),\ldots,\epsilon(d))^\transpose\in\RL^d$
is Gaussian noise with distribution $N(0,\sigma^2I)$.
The noise variance
$\sigma^2$ is assumed to be fixed and known,
and $I$ denotes the $d$-dimensional identity 
matrix. The $d\times d$ matrix $A$ is
described in terms of its {\em structure} $S$
and {\em weights} $\ww$,
where $S$ is the binary matrix
with entries $S_{ij}=\IND\{A_{ij}\neq 0\}$
and $\ww$ is the concatenation
of all the nonzero entries of $A$.
We naturally view $S$ as the {\em model}
and the weights $\ww$ as the 
associated {\em parameters}.

Throughout, we assume that $S$
is the adjacency matrix of a
directed acyclic graph (DAG) on $\{1,\ldots,d\}$,
and we often identify $S$ with the DAG it represents.
To ensure identifiability~\cite{Peters_2013},
in all our results we will assume
{\em causal minimality}~\cite{zhang:08},
which is a mild {\em faithfulness} assumption~\cite{sprites_book_2017}.
Causal minimality means that 
the observations are generated
by a model $S^*$ and parameters $\ww^*$ such that 
the law of $X$ in~(\ref{eq:LSE}) under $(S^*,\ww^*)$
is not the same as its law under $(S,\ww)$
for any subgraph $S$ of $S^*$ and any collection
of associated
parameters, $\ww$.
Note that, often,
this model is equivalently written 
as $X=A^\transpose X+\epsilon$,
instead of the form~(\ref{eq:LSE})
adopted here.

If we let $P_j(S,X)$, 
with `$P$' for `parents',
denote the sub-vector of $X$ consisting
of those of its components $X(i)$ that influence the value of
$X(j)$ in~(\ref{eq:LSE}), i.e., such that
$S_{ji}\neq 0$, then
the $j$th equation
in~(\ref{eq:LSE}) is
\begin{equation}
    X{(j)} = \mathbf{w}{(j)}^\transpose P_j(S, X) + \epsilon{(j)},  \quad j = 1, \ldots, d,
    \label{eq:model_2_by_dim}
\end{equation}
where $\ww(j)$ is the column vector that consists of
the nonzero elements of the $j$th row of $A$.
In this notation,
$\ww$ can be viewed as the concatenation of the vector
blocks $\ww(j)$, as
$\mathbf{w} 
= (\mathbf{w}{(1)}^\transpose, \ldots, 
\mathbf{w}{(d)}^\transpose)^\transpose.$

Clearly, any matrix $A$ in~(\ref{eq:LSE})
can be described by its structure $S$ and
weights $\ww$. We write 
$A=m(S,\ww)$
for this map, and we adopt the following
spike-and-slab prior on the space of matrices
$A$. First, a uniform prior
$\pi(S)\propto 1$, is used for $S$ in the space $\clG_d$
of all adjacency matrices of DAGs on $\{1,\ldots,d\}$
and then, 
given $S$, we place independent
Gaussian 
$N(0,\sigma_w^2)$ priors on the 
associated weights $\ww$.

We will find it useful to think of $AX$
as vector of linear combinations of the weights~$\ww$,
\begin{equation}
AX=M(S,X)\ww,
\label{eq:modelM}
\end{equation}
for an appropriate matrix $M(S,X)$.
Comparing~(\ref{eq:modelM}) with~(\ref{eq:LSE}) 
it is not hard to see that $M(S, X)$ is
\begin{equation*}
\begin{bmatrix}
\overbrace{P_1(S, X)^{\transpose} }^{s_1} & \!0  &\! \cdots & 0 & 0 \\
\overbrace{0 \cdots 0 }^{s_1} & \!\overbrace{  P_2(S, X)^{\transpose}  }^{s_2}  &\! \cdots & 0 & 0 \\
\vdots & \!\vdots  &\! \ddots & \vdots & \vdots \\
\overbrace{ 0 \cdots 0}^{s_1} & \!\overbrace{ 0 \cdots 0}^{s_2}  &\! \cdots & \overbrace{ 0 \cdots 0}^{s_{d-1}} & \overbrace{ P_d(S, X)^{\transpose}  }^{s_{d}}
\end{bmatrix}\!.
\end{equation*}
Here,
$s_j$ denotes the dimension of each $P_j(S,X)$, where
the vectors $P_j(S,X)$
are defined in~(\ref{eq:model_2_by_dim}).

\subsection{The posterior}

Given independent and identically
distributed (i.i.d.)\ observations $\Xn=(X_1,\ldots,X_n)$
from the model~(\ref{eq:LSE}),
and writing $\phi_v$ for the $N(0,v)$ density,
the associated likelihood is
$$f(\mathbf{X}^n|S,\mathbf{w})= \prod_{i=1}^n 
\phi_{\sigma^2}(X_i-m(S,\textbf{w})X_i).$$

First, we examine
the full conditional density of the parameters:
$\pi(\ww|S,\Xn)$ is proportional to
\begin{align*}
&
	f(\mathbf{X}^n|S,\mathbf{w})\pi(\mathbf{w}|S) \pi(S)\\
&\propto 
    \exp\Big\{  -\frac{1}{2\sigma^2}\sum_{i=1}^n \|X_i-M(S, X_i)\mathbf{w}\|_2^2  \\ &\qquad \qquad \qquad \qquad \qquad \qquad \qquad - \frac{1}{2\sigma_w^2} \|\mathbf{w}\|_2^2   \Big\} \\
&\propto 
	\exp\Big\{ -\frac{1}{2}\mathbf{w}^\transpose 
	\Big(\Sigma_M(S, \mathbf{X}^n) + \frac{1}{\sigma_w^2}I\Big)\mathbf{w}
	 \\ &\qquad \qquad \qquad \qquad \qquad \qquad \qquad  + \mathbf{w}^\transpose b(S, \mathbf{X}^n) \Big\}, 
\end{align*}
where,
\begin{align*} 
\Sigma_M(S, \mathbf{X}^n) & := \frac{1}{\sigma^2}\sum_{i=1}^n 
	M(S, X_i)^\transpose M(S, X_i),\\
 b(S, \mathbf{X}^n) &:= \frac{1}{\sigma^2}\sum_{i=1}^nM(S, X_i)^\transpose X_i .
\end{align*}
Therefore, 
$\mathbf{w}|S, \mathbf{X}^n \sim
  N(\mu_w(S,\Xn), \Sigma_{w}(S, \mathbf{X}^n)),$
where,
\begin{align*}
\Sigma_{w}(S, \mathbf{X}^n)   
&:=  
	\Big( \frac{1}{\sigma_w^2}I + \Sigma_M(S, \mathbf{X}^n) \Big)^{-1},\\
\mu_w(S,\Xn)
&:=
	\Sigma_{w}(S, \mathbf{X}^n)b(S, \mathbf{X}^n).
\end{align*}

Note that, because of the structure of $M(S,X)$,
$\Sigma_w(S,\Xn)$ is a block diagonal matrix, whose
$j$th block, $\Sigma^{(j)}_{w}(S,\Xn)$,
$1\leq j\leq d$, is given by 
\begin{equation}
	\Big(\frac{1}{\sigma_w^2}I+ \frac{1}{\sigma^2}\sum_{i=1}^n 
	P_j(S, X_i)P_j(S, X_i)^\transpose\Big)^{-1}.
	\label{eq:Sigmawdef}
\end{equation}
Similarly, for $1\leq j\leq d$,
the vector $b(S,\Xn)$ can be written in
terms of the vector blocks
\begin{equation}
b^{(j)}(S, \mathbf{X}^n) 
:= 
	\frac{1}{\sigma^2} \sum_{i=1}^n X_i(j)P_j(S, X_i),
\label{eq:bj}
\end{equation}
as $b(S,\Xn)=(b^{(1)}(S,\Xn)^\transpose,\ldots,
b^{(d)}(S,\Xn)^\transpose)^\transpose$,
and $\mu_w(S,\Xn)$ can be written in 
terms of 
\begin{equation}
\mu_w^{(j)}(S,\Xn)
:=\Sigma_w^{(j)}(S,\Xn)b^{(j)}(S,\Xn),
\label{eq:muwj}
\end{equation}
as $\mu_w(S,\Xn)=(\mu_w^{(1)}(S,\Xn)^\transpose,\ldots,
\mu_w^{(d)}(S,\Xn)^\transpose)^\transpose$.

Next, we examine the posterior on model space.
A simple computation shows that $\pi(S|\Xn)$ is
\begin{align}
&\int f(\mathbf{X}^n|S,\mathbf{w})\pi(\mathbf{w}|S)\, d\mathbf{w}
    \nonumber\\
  &\propto \exp\Big(\frac{1}{2}b(S, \mathbf{X}^n)^\transpose 
\Sigma_{w}(S, \mathbf{X}^n) b(S, \mathbf{X}^n)\Big) 
    \nonumber\\
&\quad \times
|\Sigma_{w}(S, \mathbf{X}^n)|^{1/2},
    \label{eq:final_integral_full_cond_post}
    \end{align}
where  $|B|$ denotes the determinant of matrix $B$.
Substituting into~(\ref{eq:final_integral_full_cond_post})
the expression for $M(S,X)$,
using the notation in~(\ref{eq:Sigmawdef})--(\ref{eq:bj}),
and simplifying, shows that $\pi(S|\Xn)$ is proportional to
\begin{align}
&\exp \Big(\frac{1}{2} \sum_{j=1}^d  b^{(j)}(S, \mathbf{X}^n)^{\transpose}
\Sigma_{w}^{(j)}(S, \Xn)  b^{(j)}(S, \mathbf{X}^n) \Big)
	\nonumber\\ & \quad \qquad \qquad \qquad \qquad\times\prod_{j=1}^d |\Sigma^{(j)}_{w}(S, \Xn)|^{1/2}.
    \label{eq:posterior_graphs_by_dimensions}
\end{align}
This expression will be our starting point in Section~\ref{sec:generalcase}.

\section{Posterior concentration: Binary case}
\label{sec:binarycase}

Before considering the posterior induced by the general
linear causal model in \eqref{eq:LSE}, in this section
we examine the simpler {\em binary} linear causal model,
where all the parameters in $\ww$ are equal to~1, i.e., 
 $A=S$:
\begin{equation}
    X = S X + \epsilon.
    \label{eq:model_1}
\end{equation} 
In this case, every model $S\in\clG_d$ is causally minimal,
so we have identifiability~\cite{Peters_2013}.
Also,
it is easy to obtain the 
form of the posterior on models, $\pi(S|\Xn)$, as:
\begin{equation}
    \pi(S|\mathbf{X}^n) \propto  f(\mathbf{X}^n|S)
    \propto \prod_{i=1}^n \phi_{\sigma^2}(X_i - S X_i).
    \label{eq:posterior_1}
\end{equation}

Our main result here states that, when the observations
are generated according to some matrix $S^*$, its posterior 
probability almost surely converges to~1 at an 
exponential rate, with exponent equal to~$1/2$.
Interestingly, the exponent is the same for all~$S^*$.

\begin{theorem}
\label{thm:S}
If $\{X_n\}$ are i.i.d.\ observations generated 
by the model~\eqref{eq:model_1} with $S=S^*$, 
then, as $n\to\infty$:
    $$\pi(S^*|\mathbf{X}^n) = 1-\exp\Big\{-\frac{n}{2}
+O(\sqrt{n\log \log n})\Big\} \quad \mbox{a.s.}$$
\end{theorem}

For the proof, we will need the following two lemmas.
Lemma~\ref{lem:relent},
 proved  in Appendix~\ref{appendix:BSC},
computes the Kullback-Leibler
divergence $D(P_{S'}\|P_S)$ between the laws $P_S,P_{S'}$
of $X$ in~(\ref{eq:model_1}) corresponding
to $S,S'$. 

\begin{lemma} 
\label{lem:relent}
For any $S\in\clG_d$, let $P_S$ denote the law of 
$X=(I-S)^{-1}\epsilon$ in~\eqref{eq:model_1}.
Then, for any pair $S,S'\in\clG_d$,
$$D(P_{S'}\|P_S)=\frac{1}{2}\big[\|(I-S)(I-S')^{-1}\|_F^2-d\big],$$
where 
$\|B\|_F^2:=\sum_{ij}B^2_{ij}$
is the Frobenius norm of~$B$.
\end{lemma}

The exponent~$1/2$ in the theorem
will be seen to be based on the following
simple computation. The proof of 
Lemma~\ref{lem:minD} is in Appendix~\ref{appendix:BSC}.

\begin{lemma} 
\label{lem:minD}
For any $G\in\clG_d$, we have, 
$$
\min_{S\in\clG_d, S\neq G} \|(I-S)(I-G)^{-1}\|_F^2=d+1,$$
or, equivalently,
$$
\min_{S\in\clG_d, S\neq G} D(P_{G}\|P_{S})=\frac{1}{2}.$$
\end{lemma}

We are now ready to prove the theorem.

\smallskip

\noindent
{\sc Proof of Theorem~\ref{thm:S}.}
Let $S^*\in\clG_d$ be fixed.
From the expression for the posterior in~(\ref{eq:posterior_1}) we have,
\begin{equation*}
    \begin{aligned}
        \pi(S^*|\mathbf{X}^n) &= \frac{f(\mathbf{X}^n|S^*)}{\sum_{S}f(\mathbf{X}^n|S)}
\\& = \frac{\prod_{i = 1}^n \phi_{\sigma^2}(X_i-S^*X_i)}{\sum_{S} \prod_{i = 1}^n \phi_{\sigma^2}(X_i-SX_i)},
    \end{aligned}
\end{equation*}
where the sums are over all $S\in\clG_d$. 
Using the form of the 
Gaussian density and the fact that each, we have 
$X_i = (I-S^*)^{-1}\epsilon_i$
from~(\ref{eq:model_1}), $\pi(S^*|\mathbf{X}^n)$ becomes,
$$	\frac{ \exp\big\{-\frac{1}{2\sigma^2} \sum_{i=1}^n \|\epsilon_i\|^2
	\big\}}
    	{\sum_S  \exp\big\{-\frac{1}{2\sigma^2}\sum_{i=1}^n \| 
	(I-S)(I-S^*)^{-1}\epsilon_i\|^2  \big\}   },
$$
and writing $D(S)= (I-S)(I-S^*)^{-1}$, this further becomes
$$ \frac{\exp\big\{-\frac{1}{2\sigma^2} \sum_{i=1}^n \|\epsilon_i\|^2\big\}}
        {\sum_S  \exp\big\{-\frac{1}{2\sigma^2} \sum_{i=1}^n \epsilon_i^\transpose D(S)^\transpose D(S) \epsilon_i\big\}}.
$$
Letting $U_n(S)$ denote the i.i.d.\ partial sums $U_n(S) =\sum_{i=1}^n Y_i(S)$, 
with
$$Y_i(S)= \frac{1}{2\sigma^2}\epsilon_i^\transpose [D(S)^\transpose D(S)-I]
\epsilon_i,$$
we can then express:
$$\pi(S^*|\mathbf{X}^n) = \Big[1+ \sum_{S \neq S^*} \exp\{-U_n(S)\}\Big]^{-1}.$$

For each $i$ and any $S$, 
it is easy to compute that,
\begin{equation}
\begin{aligned}
    \nu(S)&:=\BBE[Y_i(S)]=\frac{1}{2}\big[\|D(S)\|_F^2 - d\big]
\\&=\frac{1}{2}\Big[\|(I-S)(I-S^*)^{-1}\|_F^2-d\Big],
\end{aligned}
\label{eq:muS}
\end{equation}
so that,
$$\pi(S^*|\mathbf{X}^n) = \Big[1+ \sum_{S \neq S^*} 
e^{-n\nu(S)}\exp\{-\bar{U}_n(S)\}\Big]^{-1},$$
where $\bar{U}_n(S)$ are the centered partial sums
$\bar{U}_n(S)=U_n(S)-n\nu(S)$.
In view of the law of the iterated logarithm (LIL), this already contains all the information
we need -- the remainder of the proof is simple asymptotic estimates.

From~(\ref{eq:muS}) and Lemma~\ref{lem:minD}, we know that the
minimum of $\nu(S)$ over all $S\neq S^*$ is  $1/2$.
Therefore,
\begin{equation}
\pi(S^*|\mathbf{X}^n) \geq \Big[1+ e^{-\frac{n}{2}}\sum_{S \neq S^*} \exp\{-\bar{U}_n(S)\}\Big]^{-1}.
\label{eq:bound1}
\end{equation}
And by the LIL we know that, for each $S\neq S^*$,
\begin{equation}
|\bar{U}_n(S)|=O(\sqrt{n\log\log n})\quad \mbox{a.s.},
\label{eq:bound2}
\end{equation}
as $n\to\infty$, since the variance,
$${\rm Var}(Y_i(S))=
{\rm Var}\Big(\frac{1}{2\sigma^2}
\epsilon_i^\transpose (D(S)^\transpose D(S) -  I)\epsilon_i\Big),$$ 
is clearly finite. 
Combining~(\ref{eq:bound1}) with~(\ref{eq:bound2}) and
the elementary bound
$(1+x)^{-1}\geq 1-x$, $x\geq 0$, 
gives
$$\pi(S^*|\mathbf{X}^n) \geq 1-\exp\big\{-\frac{n}{2}+
O(\sqrt{n\log \log n})\big\} \quad \mbox{a.s.}$$

To obtain a matching upper bound, we  note that
$$\pi(S^*|\Xn) \leq \big[1+ e^{-n/2} \exp\{-\bar{U}_n(S_1)\}\big]^{-1}, $$
for any $S_1\in \mathcal{G}_d$, $S_1\neq S^*$,
and hence,
$$\pi(S^*|\mathbf{X}^n) \leq 1-\exp\big\{-\frac{n}{2}
+O(\sqrt{n\log \log n})\big\} \quad \mbox{a.s.}$$
where we used 
the elementary bound $(1+x)^{-1} \leq 1-x+x^2$, $x \geq 0$.
\qed

\section{Posterior concentration: General case}
\label{sec:generalcase}

In this section we consider the posterior $\pi(S|\Xn)$ 
on models, under the general linear causal model given
in~(\ref{eq:LSE}) and ~(\ref{eq:model_2_by_dim}). 
We assume throughout that $\{X_n\}$ are i.i.d.\
observations from the model $X=A^*X+\epsilon$,
with respect to the given true underlying matrix
$A^*$ with structure $S^*\in\clG_d$
and associated parameters~$\ww^*$.
Theorems~\ref{thm:final}
and~\ref{thm:final2} contain the main results of this paper.
They show that, if $S^*$ is
maximal, then $\pi(S^*|\Xn)$ converges to 1 exponentially
fast, otherwise, it converges at 
rate no faster than $1/\sqrt{n}$.

Recall the expression for the 
posterior $\pi(S|\Xn)$ given in~(\ref{eq:posterior_graphs_by_dimensions}).
Writing $T^{(j)}(S,\Xn)$,
namely, $1\leq j\leq d$,
for the $j$th summand in the 
exponent 
in~(\ref{eq:posterior_graphs_by_dimensions}),
namely, 
\begin{equation}
b^{(j)}(S, \mathbf{X}^n)^\transpose\Sigma_w^{(j)}
(S, \mathbf{X}^n)   b^{(j)}(S, \mathbf{X}^n),
\label{eq:Tj}
\end{equation}
we see that $\pi(S|\Xn)$ is proportional to: 
\begin{equation}
    \exp \Bigg( \frac{1}{2} \sum_{j=1}^d 
T^{(j)}(S, \Xn) \Bigg)\prod_{j=1}^d \sqrt{ |\Sigma_{w}^{(j)}(S, \Xn)|}.
\label{eq:posterior_new}
\end{equation}

We first derive the first-order asymptotic behavior 
of 
$\Sigma_w^{(j)}(S,\Xn)$ and
$T^{(j)}(S, \Xn)$.
Proposition~\ref{prop:exponent_approx_A} is proved in
Appendix~\ref{app:GLS}.

\begin{proposition}
\label{prop:exponent_approx_A}
Suppose $\{X_n\}$ are i.i.d.\  observations
generated the model~\eqref{eq:LSE}, 
with respect to a matrix $A^*$ with 
model $S^*$
and parameters $\mathbf{w}^*$,
such that $(S^*,\ww^*)$ is causally minimal.
Then, for any $S\in\clG_d$ and each $1\leq j\leq d$, 
as $n\to\infty$, we have,
almost surely,
\begin{align*}
n\Sigma_{w}^{(j)}(S,\Xn) 
&\to 
	\Sigma_{\infty}^{(j)}(S),\\
T^{(j)}(S, \Xn) 
&=  
	n T_\infty^{(j)}(S) 
	+ O(\sqrt{n\log \log n}),
\end{align*}
where,
$$\Sigma_{\infty}^{(j)}(S) 
:= 
	\sigma^2\mathbb{E}[P_j(S, X)P_j(S, X)^\transpose]^{-1},
$$
and $T_{\infty}^{(j)}(S)$ is given by
\begin{align*}
&
	\sigma^2\big(\mathbb{E}[X{(j)}P_j(S, X)^\transpose]\\
&  \qquad\mathbb{E}
	[P_j(S, X)P_j(S, X)^\transpose]^{-1}\mathbb{E}[X{(j)}P_j(S, X)] \big),
\end{align*}
or, equivalently,
$$
	\mathbb{E}[X{(j)^2}] - \min_{\eta \in \mathbb{R}^{s_j}}
	\mathbb{E}\Big[\big(X{(j)}- \eta^\transpose P_j(S, X)\big)^2\Big].
$$
\end{proposition}

Next, recall from~(\ref{eq:muwj}) the definition 
of the $j$th component, $\mu_w^{(j)}(S,\Xn)$, of the mean
of $\ww$ under the the full conditional 
distribution $\pi(\ww|S, \Xn)$. 
Proposition~\ref{prop:exponent_approx} is
proved in Appendix~\ref{app:GLS}.

\begin{proposition}
\label{prop:exponent_approx}
Under the assumptions of Proposition~\ref{prop:exponent_approx_A},
for any $S\in\clG_d$ and each $1\leq j\leq d$, as $n\to\infty$, we have,
almost surely,
$$\mu_w^{(j)}(S, \Xn) 
\to
	\mu^{(j)}_\infty(S),
$$
where, $\mu_\infty^{(j)}(S)$ is given by:
\begin{align}
&
	\mathbb{E}[P_j(S, X)P_j(S, X)^\transpose]^{-1} 
	\mathbb{E}[X{(j)}P_j(S, X)]
	\nonumber\\
&=
	\argmin_{\eta \in \mathbb{R}^{s_j}}\mathbb{E}\Big[
	\big(X{(j)}- 
	\eta^\transpose P_j(S, X)\big)^2\Big].
    \label{eq:def_T_i}
\end{align}
In particular,
under the true model $S=S^*$, 
we have $\mu^{(j)}_\infty(S^*) = \mathbf{w}^*{(j)}$,
$j=1,\ldots,d$.
\end{proposition}

As with $\mu(S, \Xn)$ and $\Sigma_w(S, \Xn)$,
we let
$\mu_\infty(S) := (\mu^{(1)}_\infty(S)^\transpose, \ldots,
\mu^{(d)}_\infty(S)^\transpose)^\transpose$,
and we write $\Sigma_\infty(S)$ for the block
diagonal matrix with blocks
$\Sigma_\infty^{(j)}(S)$, $1\leq j\leq d$.
Also, for any matrix $A$ in~(\ref{eq:LSE}),
let $P_A$ denote the law of $X=(A-I)^{-1}\epsilon$.
The following is a straightforward generalization
of Lemma~\ref{lem:relent}. Its proof is identical to that
of Lemma~\ref{lem:relent}, and hence omitted.

\begin{lemma} 
\label{lem:relent2}
For any pair of matrices $A,A'$ with models
$S,S'\in\clG_d$, respectively, we have:
$$D(P_{A'}\|P_A)=\frac{1}{2}\big[\|(I-A)(I-A')^{-1}\|_F^2-d\big].$$
\end{lemma}

Recall that we call $S^*$ {\em maximal} if
adding any one more edge to it would necessarily violate
the acyclicity constraint.
In the following two lemmas, we estimate the posterior
log-odds between $S^*$ and different models $S\in\clG_d$. 
These results, 
proved in Appendix~\ref{app:mainlemmas}, are the main technical
ingredients in the proofs of our two main results,
Theorems~\ref{thm:final}
and~\ref{thm:final2} below.

\begin{lemma}
\label{lemma:maximal_case}
Under the assumptions of Proposition~\ref{prop:exponent_approx}, 
for any $S\in\clG_d$ such that $S^*$ is not a subgraph of $S$,
as $n\to\infty$ we have, almost surely,
\begin{align*}
\log\frac{\pi(S|\mathbf{X}^n)}{\pi(S^*|\mathbf{X}^n)} 
&= 
-nD(P_{A^*} \| P_{A(S)})\\
&\quad+O(\sqrt{n\log \log n})
\end{align*}
where $A(S) = m(S, \mu_\infty(S))$ and $D(P_{A^*}\|P_{A(S)})>0$.
\end{lemma}

\begin{lemma}
\label{lemma:non_maximal_case}
Under the assumptions of Proposition~\ref{prop:exponent_approx},
and assuming $S^*$ is non-maximal, there is an
$S^+ \in \clG_d$ such that $S^*$ is a subgraph of $S^+$
and $S^*$ differs from $S^+$ in only one edge. 
For that $S^+$ we have,
\begin{equation*}
    \log\frac{\pi(S^+|\mathbf{X}^n)}{\pi(S^*|\mathbf{X}^n)}  
= -\frac{\log n}{2} +\delta_n +o(1) \quad \mbox{a.s.},
\end{equation*}
as $n\to\infty$,
where $\{\delta_n\}$ is eventually a.s.\ nonnegative.
\end{lemma}

We are now ready to state our main results.

\begin{theorem}[Posterior concentration for maximal models]
\label{thm:final}
Under the assumptions of Proposition~\ref{prop:exponent_approx},
if the true underlying model $S^*$ is maximal,
then, 
$$\pi(S^*|\mathbf{X}^n)= 
1-\exp\Big\{-nD(A^*)+
O(\sqrt{n\log \log n})\Big\},$$
almost surely as $n\to\infty$, where $D(A^*)>0$ is
\begin{equation}
D(A^*):=\min_{S \in \mathcal{G}_d: S\neq S^*}
D(P_{A^*}\|P_{m(S, \mu_\infty(S))}).
\label{eq:minS}
\end{equation}
\end{theorem}

\noindent
{\sc Proof.}
Since $S^*$ is maximal, it is not a subgraph of any 
$S\in\clG_d$, $S\neq S^*$.
Therefore, by Lemma~\ref{lemma:maximal_case}, for any 
$S\in\clG_d$, $S\neq S^*$, we have, almosto surely as $n\to\infty$,
\begin{align}
\log\frac{\pi(S|\mathbf{X}^n)}{\pi(S^*|\mathbf{X}^n)} &= 
-nD(P_{A^*}\| P_{A(S)} ) \\
&\quad+O\big(\sqrt{n\log \log n}\big).
\label{eq:maximalpre}
\end{align}
Let $S^-$ achieve the minimum in~(\ref{eq:minS}).
Then,
$$\pi(S^*|\Xn)  
= 
\frac{1}{1+\sum_{S\neq S^*}\frac{\pi(S|\Xn)}{\pi(S^*|\Xn)}}
\leq
\frac{1}{1+\frac{\pi(S^-|\Xn)}{\pi(S^*|\Xn)}},$$
and the elementary bound $(1+x)^{-1}\leq1-x+x^2$,
$x\geq 0$, combined with~(\ref{eq:maximalpre}) gives,
a.s.\ as $n\to\infty$,
$$\pi(S^*|\Xn)  
\leq
1-
\exp\Big\{-nD(A^*)+
O(\sqrt{n\log \log n})\Big\}.
$$
Similarly, $\pi(S^*|\Xn)$ is bounded below by
\begin{align*}
&	\Big[1+|\clG_d|\max_{S\neq S^*}\frac{\pi(S|\Xn)}{\pi(S^*|\Xn)}\Big]^{-1}\\
& = 
	\Big[
	1+|\clG_d|\exp
	\Big\{
	-n\min_{S\neq S^*}D(P_{A^*}\| P_{A(S)})	\\
	&\qquad \qquad \qquad\qquad \qquad \quad +O(\sqrt{n\log \log n})
	\Big \}\Big]^{-1} \\
& \geq 1-\exp\Big\{-nD(A^*)+O(\sqrt{n\log \log n})\Big\} \quad \mbox{a.s.},
\end{align*}
where we used~(\ref{eq:maximalpre}) and the simple bound
$(1+x)^{-1}\geq 1-x$, $x\geq 0$.
\qed

\begin{theorem}[Posterior concentration for non-maximal models]
\label{thm:final2}
Under the assumptions of Proposition~\ref{prop:exponent_approx},
if the true underlying model $S^*$ is non-maximal, then
its posterior probability, $\pi(S^*|\Xn)$, converges
to~1 a.s.\ as $n\to\infty$, but this convergence
is no faster than at rate $1/\sqrt{n}$.
Specifically, as $n\to\infty$ we have, 
$$1-\pi(S^*|\mathbf{X}^n) \geq (1+o(1))\frac{1}{\sqrt{n}}
\quad \mbox{a.s.}$$
\end{theorem}

\noindent
{\sc Proof.} 
Consistency follows from standard Bayesian model-selection 
theory \`{a} la Schwartz~\cite{schwartz:65}.
For more specific (and more general) results in the present 
context see, e.g.,~\cite{chickering2002optimal,cao:19}.
For the bound on the rate, 
with $S^+$  as in Lemma~\ref{lemma:non_maximal_case},
note that
$$\pi(S^*|\Xn) 
= 
	\frac{1}{1+\sum_{S\neq S^*}\frac{\pi(S|\Xn)}{\pi(S^*|\Xn)}}\\
\leq 
	\frac{1}{1+\frac{\pi(S^+|\Xn)}{\pi(S^*|\Xn)}}.
$$
Then, using
Lemma~\ref{lemma:non_maximal_case} together with the
elementary bound $(1+x)^{-1}\leq 1-x+x^2$, $x\geq 0$,
yields,
$$\pi(S^*|\Xn)\leq 1-\frac{1}{\sqrt{n}}(1+o(1))+O\Big(\frac{1}{n}\Big),
$$
as required.
\qed

\section{Edge detection}
\label{sec:edge_detection}

In this section, we draw an interesting and potentially
practically useful connection between the Bayesian
causal discovery framework of Section~\ref{sec:intro}
and optimal hypothesis testing for edge detection.
Given $n$ i.i.d.\ observations $\mathbf{X}^n=(X_1,\ldots,X_n)$ from 
the general linear causal model~(\ref{eq:LSE}),
the goal here is to estimate the support $\bchi(S) = S + S^{\transpose}$
 of the true underlying model~$S$. 

Suppose that, after $S$ is drawn
according to the prior $\pi(S)$ on $\clG_d$,
the observations $\Xn$ are (conditionally)
i.i.d.\ samples 
from~(\ref{eq:LSE}).
Given a {\em detector} 
$\hbchi=\hbchi(\Xn)$,
following~\cite{NPCausalISIT} we define the \textit{false positive} 
and \textit{false negative} 
detection rates, respectively, as,
\begin{align}
\varepsilon_n^+ &= \frac{\mathbb{E}\big[\sum_{i>j}\mathbbm{I}\{\hbchi_{ij}=1, \bchi_{ij}=0\}\big]}{\mathbb{E}\big[\sum_{i>j}\mathbbm{I}\{\bchi_{ij}=0\}\big]} 
\label{def::falsepositive} ,
    \\
\varepsilon_n^- &= \frac{\mathbb{E}\Big[\sum_{i>j}\mathbbm{I}\{\hbchi_{ij}=0, \bchi_{ij}=1\}\Big]}{\mathbb{E}\Big[\sum_{i>j}\mathbbm{I}\{\bchi_{ij}=1\}\Big]}\label{def::falsenegative}.
\end{align}
We aim to minimize $\varepsilon_n^-$, subject to $\varepsilon_n^+$ being 
below a given 
pre-specified user tolerance
$\alpha$:
\begin{equation}\label{def::NPCI}
\begin{aligned}
\inf_{\hat{\boldsymbol{\chi}}} \varepsilon_n^- \; \; \; \; 
\textrm{s.t.} \;\; \varepsilon_n^+ \leq \alpha.
\end{aligned}
\end{equation}

Let $\clC_{ij}$ denote the set of all $S\in\clG_d$ such 
that $\bchi_{ij}(S)=[S+S^\transpose]_{ij}=0$,
and let $\rho(\cdot|S)$ denote the marginal likelihood
of $\Xn$ given model $S$.
Then $\rho(\cdot|S)$ has density:
$$
    f(\Xn|S) = \int f(\Xn|S,\mathbf{w})\pi(\mathbf{w}|S) d\mathbf{w}.
$$
We also define the probability measures,
\begin{align*}
\bar{\clP}_{ij,n}(\cdot)
=
    \sum_{S\in\clC_{ij}}\frac{\rho(\cdot|S)}{\pi(\clC_{ij})} \pi(S),
\\
\bar{\clQ}_{ij,n}(\cdot)
=
    \sum_{S\in\clC_{ij}^c}\frac{\rho(\cdot|S)}{\pi(\clC_{ij}^c)} \pi(S),
 \end{align*}
that describe the law of $\Xn$ conditional on $S\in\clC_{ij}$ 
and on $S\not\in\clC_{ij},$ respectively.
Let,
$$    u^+_{ij} = \frac{\pi(\clC_{ij})}{\sum_{i>j}\pi(\clC_{ij})}, \qquad
    u^-_{ij} = \frac{\pi(\clC_{ij}^c)}{\sum_{i>j}\pi(\clC_{ij}^c)}.
$$
Then for any detector
$\hbchi=\hbchi(\Xn)$
we have:
\begin{align*}
       \varepsilon_n^+ &= \sum_{i>j} u^+_{ij}\bar{\clP}_{ij,n}(\hbchi_{ij} = 1),
        \\
       \varepsilon_n^- &= \sum_{i>j} u^-_{ij}\bar{\clQ}_{ij,n}(\hbchi_{ij} = 0).
\end{align*}

Letting $p_{ij,n}$ and $q_{ij,n}$ denote the densities 
of $\bar{\clP}_{ij,n}$ and $\bar{\clQ}_{ij,n}$, respectively,
the optimal detector, identified in~\cite{NPCausalISIT},  is given by
\begin{equation}
    \hbchi_{ij}= 0 \iff \frac{p_{ij,n}(\Xn)}{q_{ij,n}(\Xn)} \geq \frac{u_{ij}^-}{u_{ij}^+} \gamma,
\label{eq:optimaldet}
\end{equation}
where $\gamma$ is chosen so that the constraint on the 
false positive rate in~\eqref{def::NPCI} is satisfied with equality. 

As discussed in~\cite{NPCausalISIT}, the optimal detector is typically computationally impractical, since the evaluation of $p_{ij,n}$ and 
$q_{ij,n}$ requires averaging over all models $S$ in $\clC$ and in $\clC^c$, 
respectively, of which there are super-exponentially many~\cite{beingBayesian}.
However, the following form for the optimal detector is more amenable 
to numerical computations, for example using MCMC
model averaging. The utility of the
representation of the optimal detector in Proposition~\ref{thm:OptDetector}
is illustrated through simulation experiments in Section~\ref{s:edgedet}.

\begin{proposition}
\label{thm:OptDetector}
The optimal detector in~\eqref{eq:optimaldet} can be equivalently expressed 
in terms of the posterior as,
    \begin{equation}
        \hbchi_{ij}= 0 \iff 
	\pi(\bchi_{ij}(S) = 0|\Xn) \geq \frac{\gamma'}{1+\gamma'},
    \label{eq:test}
    \end{equation}
    where $\gamma' = \frac{u_{ij}^-\pi(\clC_{ij})}{u_{ij}^+\pi(\clC_{ij}^c)} \gamma$.
\end{proposition}

\noindent
{\sc Proof.}
Since,
$f(\Xn|S) = \pi(S|\Xn) f(\Xn),$
where $f(\Xn)=\sum_S\int f(\Xn|S, \mathbf{w})\pi(\ww|S)\pi(S)\,d\ww$
is the marginal likelihood, we have
\begin{align*}
    p_{ij,n}(\Xn)
&=
    \frac{f(\Xn)}{\pi(\clC_{ij})}\sum_{S\in\clC_{ij}} \pi(S|\Xn)
    \\
    &= \frac{f(\Xn)}{\pi(\clC_{ij})}\sum_{S}  \mathbbm{I}\{\bchi_{ij}(S) = 0\}\pi(S|\Xn) 
    \\
    &= \frac{f(\Xn)}{\pi(\clC_{ij})}\BBE[\mathbbm{I}\{\bchi_{ij}(S) = 0\}|\Xn].
\end{align*}
Similarly, 
$$    q_{ij,n}(\Xn) = \frac{f(\Xn)}{\pi(\clC_{ij}^c)}
	(1-\BBE[\mathbbm{I}\{\bchi_{ij}(S) = 0\}|\Xn]).
$$
Therefore, we have $\hbchi_{ij}= 0$ if and only if,
$$\frac{\BBE[\mathbbm{I}\{\bchi_{ij}(S) = 0\}|\Xn]}{1-\BBE[\mathbbm{I}\{\bchi_{ij}(S) = 0\}|\Xn]} \geq \frac{u_{ij}^-\pi(\clC_{ij})}{u_{ij}^+\pi(\clC_{ij}^c)} \gamma.
$$
The definition of $\gamma'$ and rearranging give~(\ref{eq:test}).
\qed

\section{Empirical results}
\label{sec:empirical_results}

In this section, we present numerical results that illustrate the accuracy of the theoretical 
results of Sections~\ref{sec:binarycase}--\ref{sec:edge_detection}.

\subsection{Posterior concentration for the binary model}

We consider 10000 observations $\{X_n\}$ generated by 
the binary causal model~(\ref{eq:model_1})
with respect to each
one of the 25 possible DAGs on three nodes.
\begin{figure*}[!ht]
\includegraphics[width=2in]{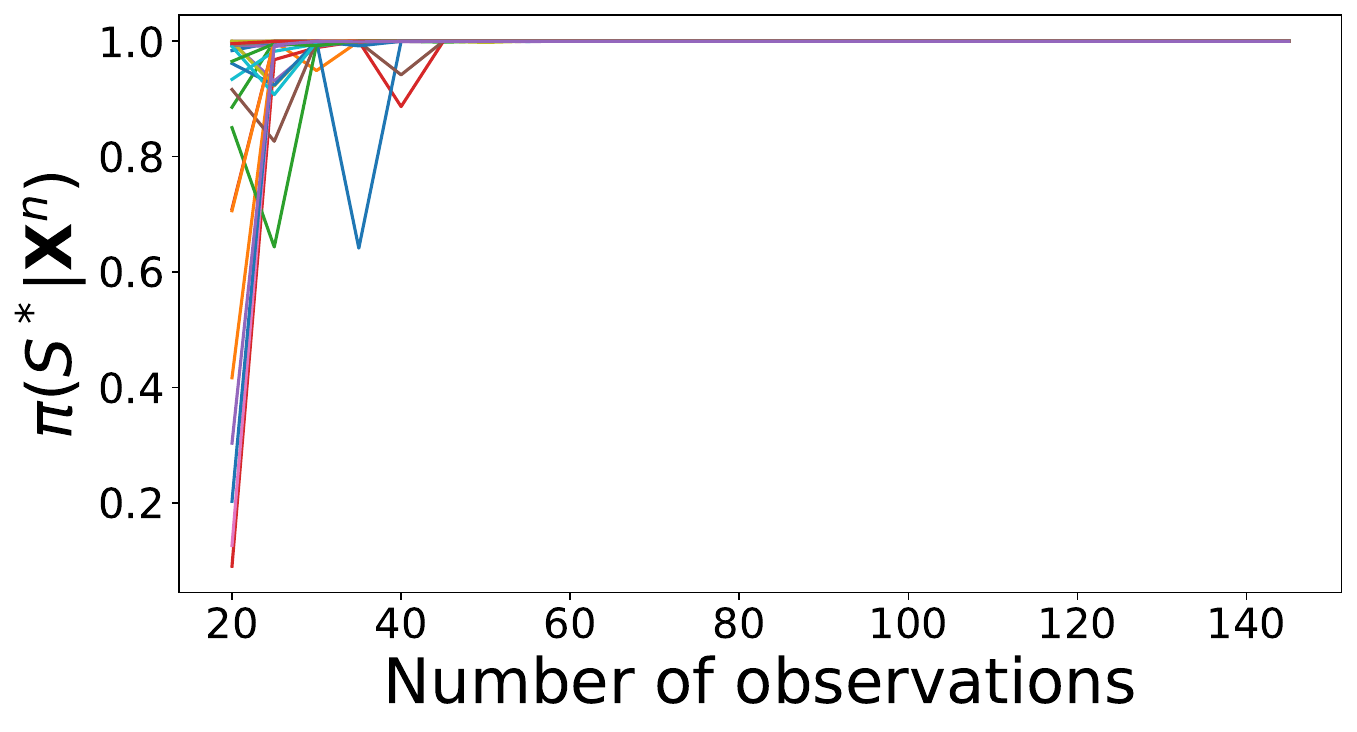}
\includegraphics[width=2.13in]{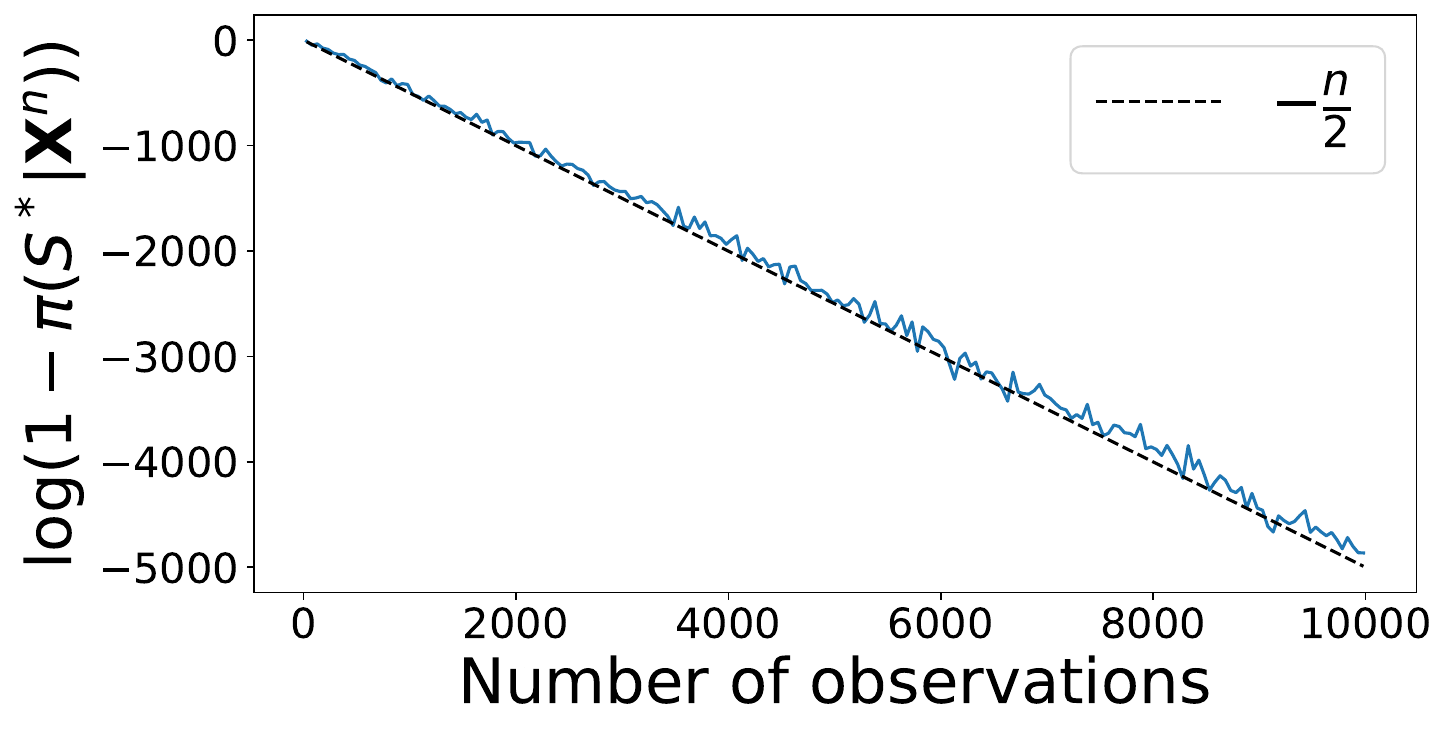}
\includegraphics[width=2.08in]{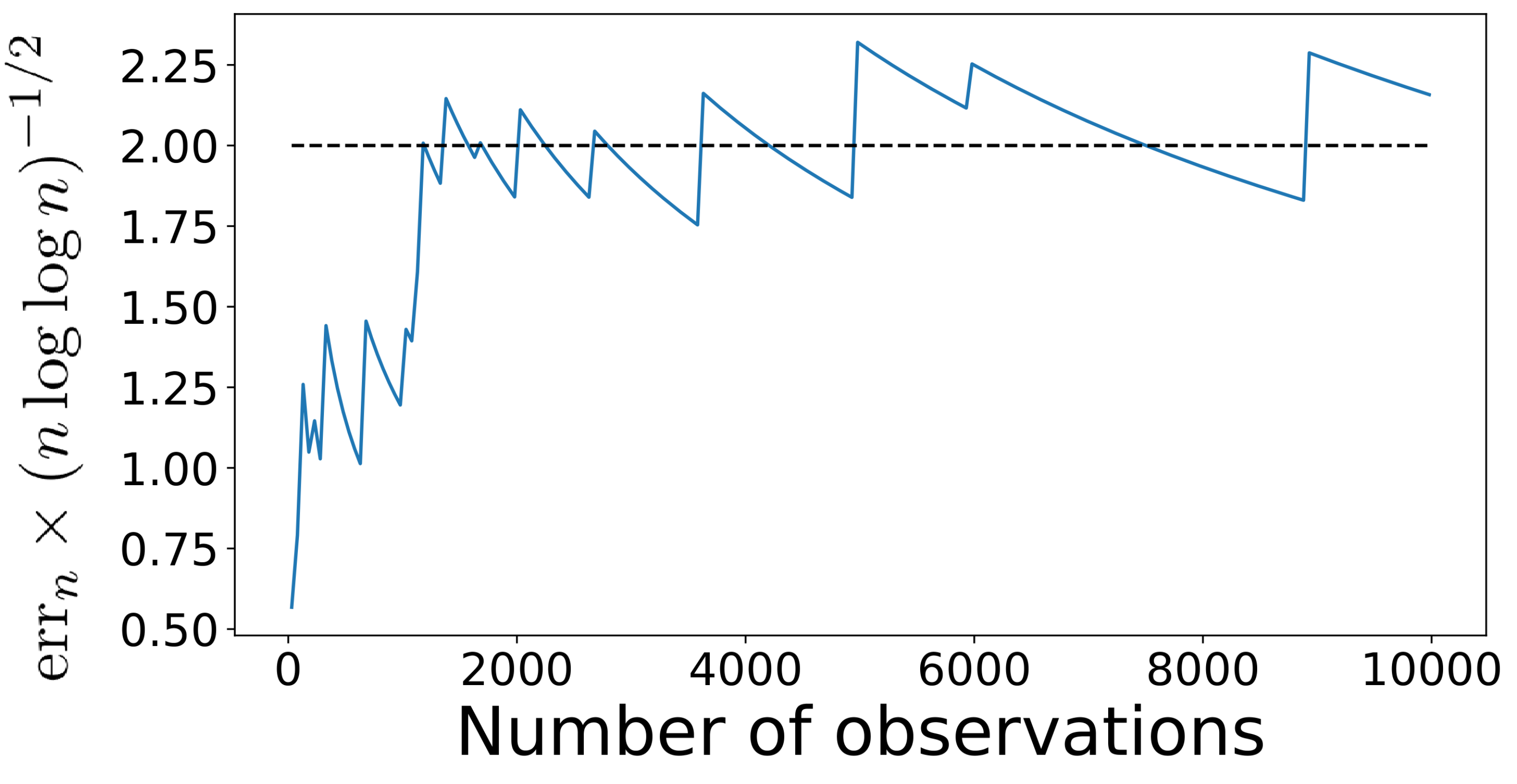}
\caption{Binary structure learning for three-node graphs. 
Left: The posterior probability $\pi(S^*|\Xn)$ of the true model $S^*$, 
for date generated by all the 25 different possible DAGs $S^*$; 
each color represents 
observations from a different DAG. 
Middle: The exponential rate of convergence of $\pi(S^*|\Xn)$ to 1.  
Right: Second order fluctuations.}
    \label{fig:binary_case}
\end{figure*}
The left plot in Figure~\ref{fig:binary_case} shows
the convergence of the unnormalized posterior
probability $\pi(S^*|\Xn)$ 
of the true model $S^*$
to~1, as a function 
of the sample size $n$, 
for each $S^*\in\clG_3$.
The middle plot 
shows $\log(1-\pi(S^*|\Xn))$ as a function on $n$, 
illustrating the first term in the 
convergence of $\pi(S^*|\Xn)$ to 1,
in the case of data generated 
by a specific $S^*$. It 
clearly indicates that the convergence is indeed
exponential with exponent $-1/2$,
as expected from Theorem~\ref{thm:S}. In the right
plot in Figure~\ref{fig:binary_case},
we examine the accuracy of the 
$O(\sqrt{n\log\log n})$ error term in Theorem~\ref{thm:S},
with simulated data from the same model.
We let,
$${\rm diff}_n:=\log(1-\pi(S^*|\Xn)+n/2,$$
and consider the fluctuations of a finite-$n$ 
approximation to the limsup
of ${\rm diff}_n$, namely,
\begin{equation}
\mbox{err}_n=\max_{1\leq k\leq n} {\rm diff_n}.
\label{eq:errn}
\end{equation}
The graph shows
${\rm err_n}\times (n\log\log n)^{-1/2}$
as a function of $n$. The resulting plot 
appears to fluctuate around an asymptotically
constant value,
suggesting that the 
$O(\sqrt{n\log\log n})$ term in Theorem~\ref{thm:S}
is of optimal order.

In all these  experiments,  the values of the posterior
$\pi(S^*|\Xn)$ are computed exactly. This is possible
because there are only 25 DAGs $S\in\clG_3$ in the
case of only three variables.

\subsection{Posterior concentration for the general model}

Here we first consider 50000 observations $\{X_n\}$
generated from the general linear model~(\ref{eq:LSE})
with respect to a matrix $A_1^*$ corresponding to a
{\em maximal} model $S_1^*$:
\begin{equation}
    A_1^* = \begin{bmatrix}
0 & 1.77 & -0.35 \\
0 & 0 & 0 \\
0 & 0.26 & 0
\end{bmatrix}\!,
\;  
S_1^* = \begin{bmatrix}
0 & 1 & 1 \\
0 & 0 & 0 \\
0 & 1 & 0
\end{bmatrix}.
\label{eq:matexample1}
\end{equation}

\begin{figure*}[!ht]
 \includegraphics[width=0.32\linewidth]{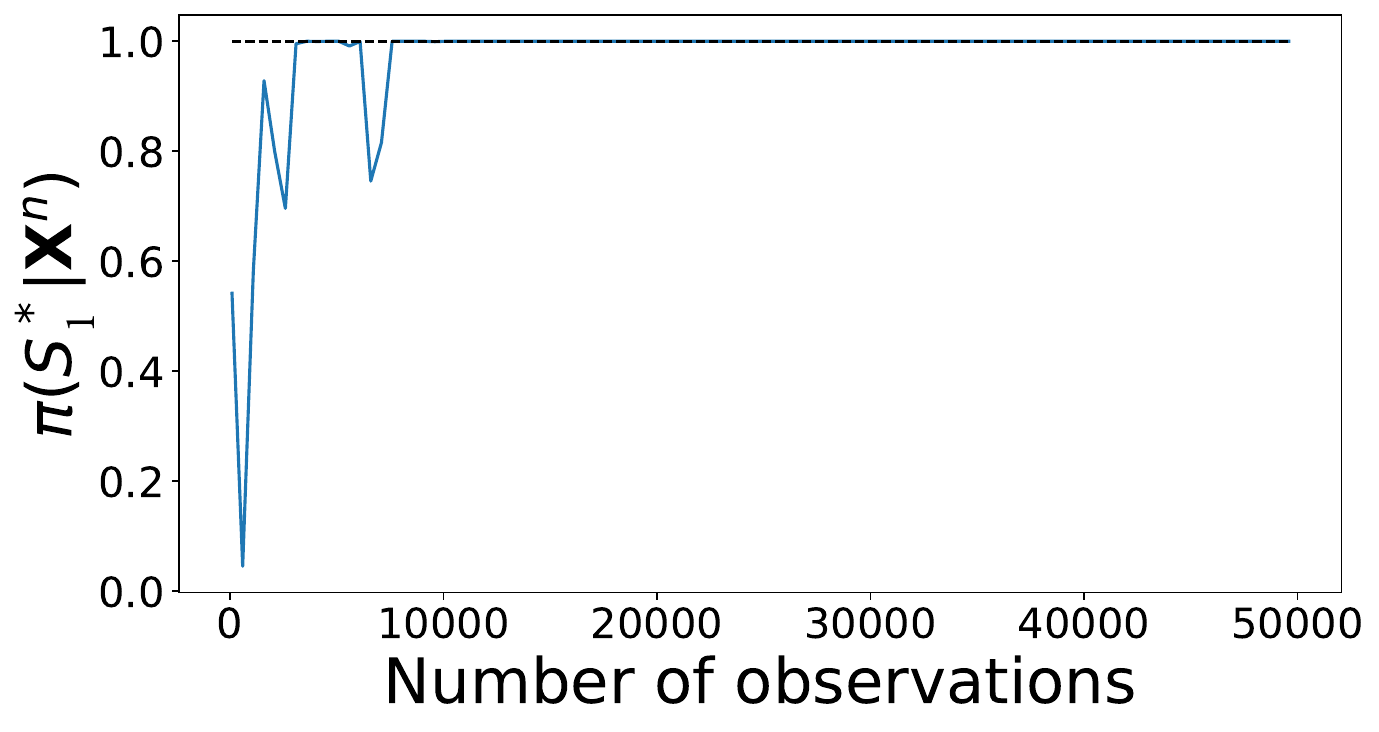}
 \includegraphics[width=0.33\linewidth]{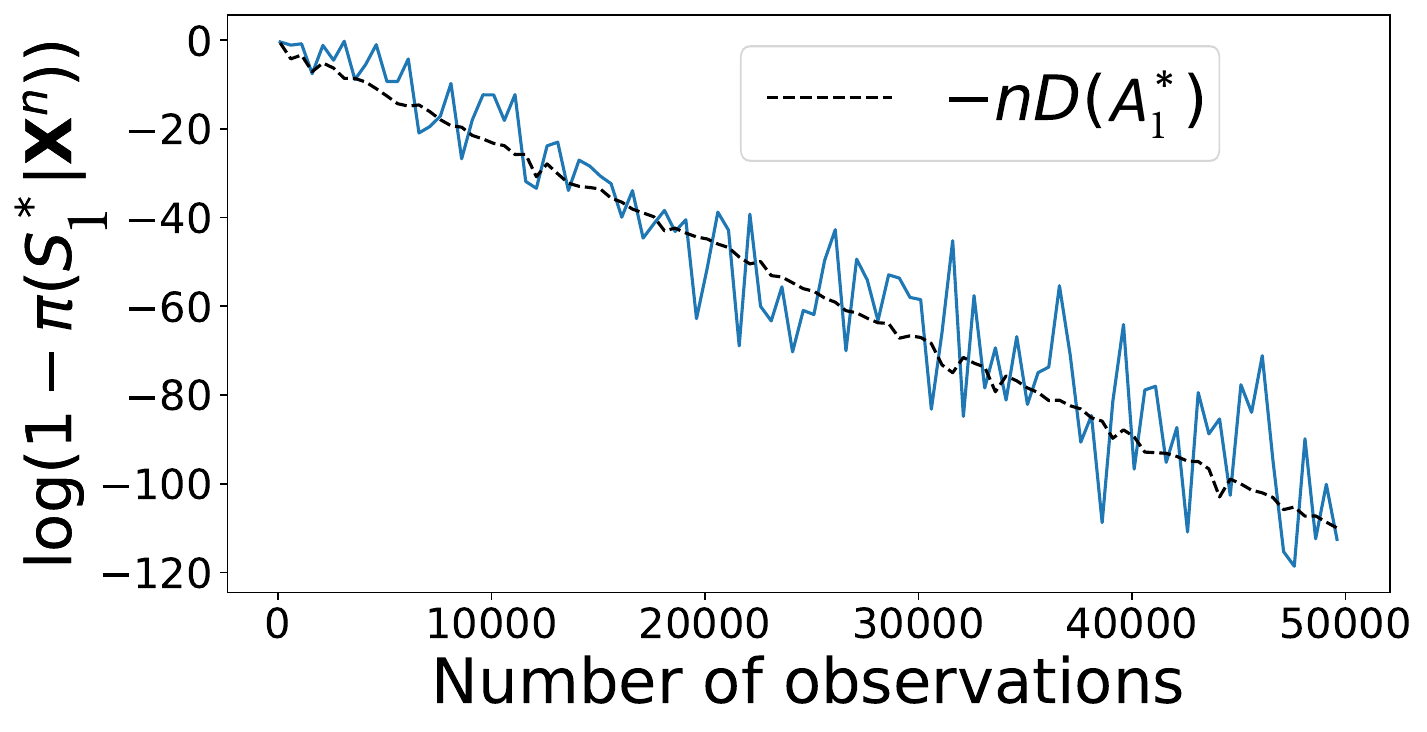}
 \includegraphics[width=0.328\linewidth]{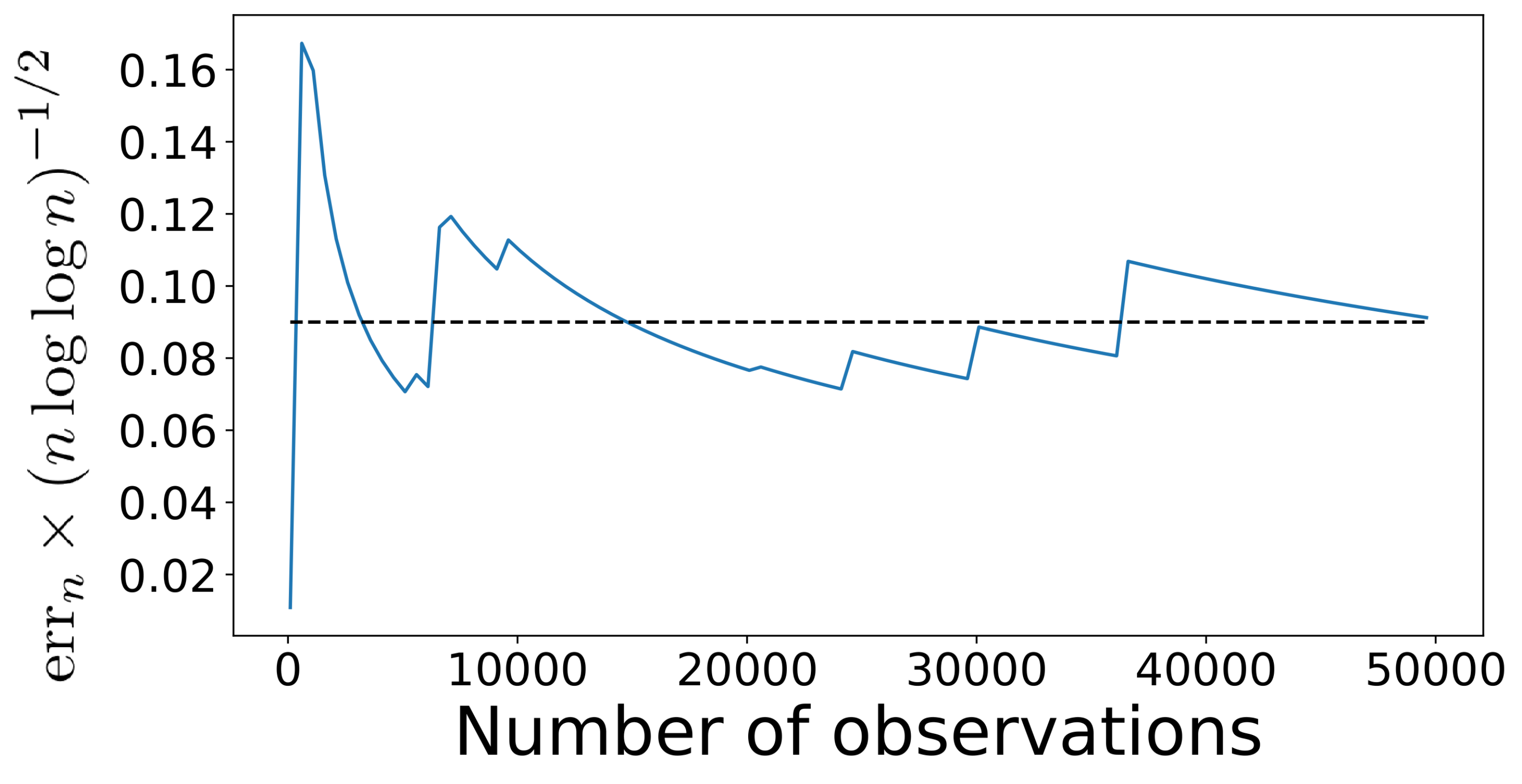}
    \caption{General structure learning for a maximal three-node graph. 
The true matrix is $A_1^*$ is given in~(\ref{eq:matexample1}).
Left: The posterior probability $\pi(S_1^*|\Xn)$ of the true model $S_1^*$.
Middle: The exponential rate of convergence of $\pi(S_1^*|\Xn)$ to 1. 
Right: Second order fluctuations.}
\label{fig:general_maximal_case}
\end{figure*}

Again, we look at the behavior of the posterior 
probability $\pi(S_1^*|\Xn)$ at three different
scales. The left plot in Figure~\ref{fig:general_maximal_case}
shows the convergence of the unnormalized probability
$\pi(S_1^*|\Xn)$ to~1, as a function of the number
of $n$ of observations between
0 and 50000. The middle plot illustrates
the first term in the rate of convergence of $\pi(S_1^*|\Xn)$ to 1,
showing $\log(1-\pi(S^*|\Xn))$ as a function of~$n$.
As expected from Theorem~\ref{thm:final}, the resulting plot
is linear with slope $-D(A_1^*)$.
In the right plot we look at the approximation
error in the exponent,
$${\rm diff}_n:=\log(1-\pi(S^*|\Xn)+D(A_1^*),$$
and consider the fluctuations of a finite-$n$ 
approximation to its limsup, given by
$\mbox{err}_n$ as in~(\ref{eq:errn}).

As before, the results indicate that the error term
$O(\sqrt{n\log\log n})$ in Theorem~\ref{thm:final}
is of optimal order.

\begin{figure*}[!ht]
\centering
\includegraphics[width=0.334\linewidth]{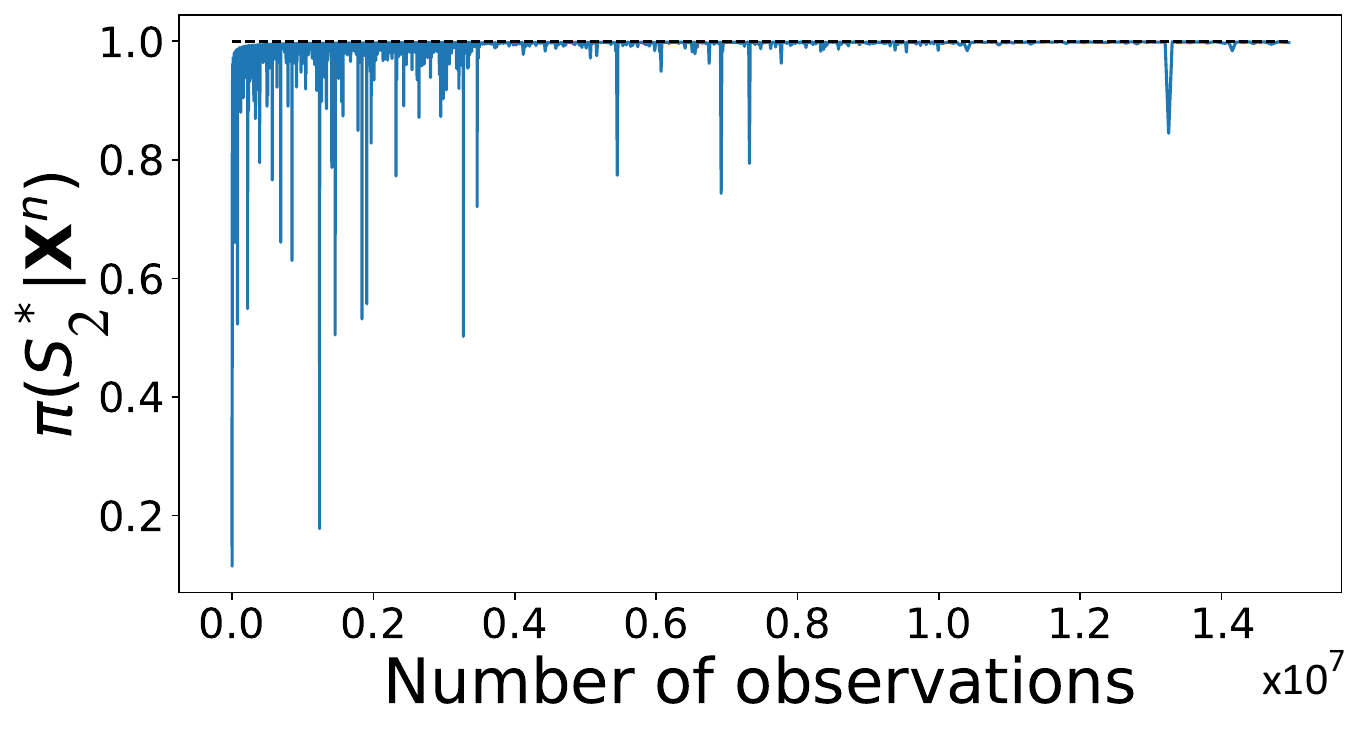}
\hspace{0.1in}
\includegraphics[width=0.33\linewidth]{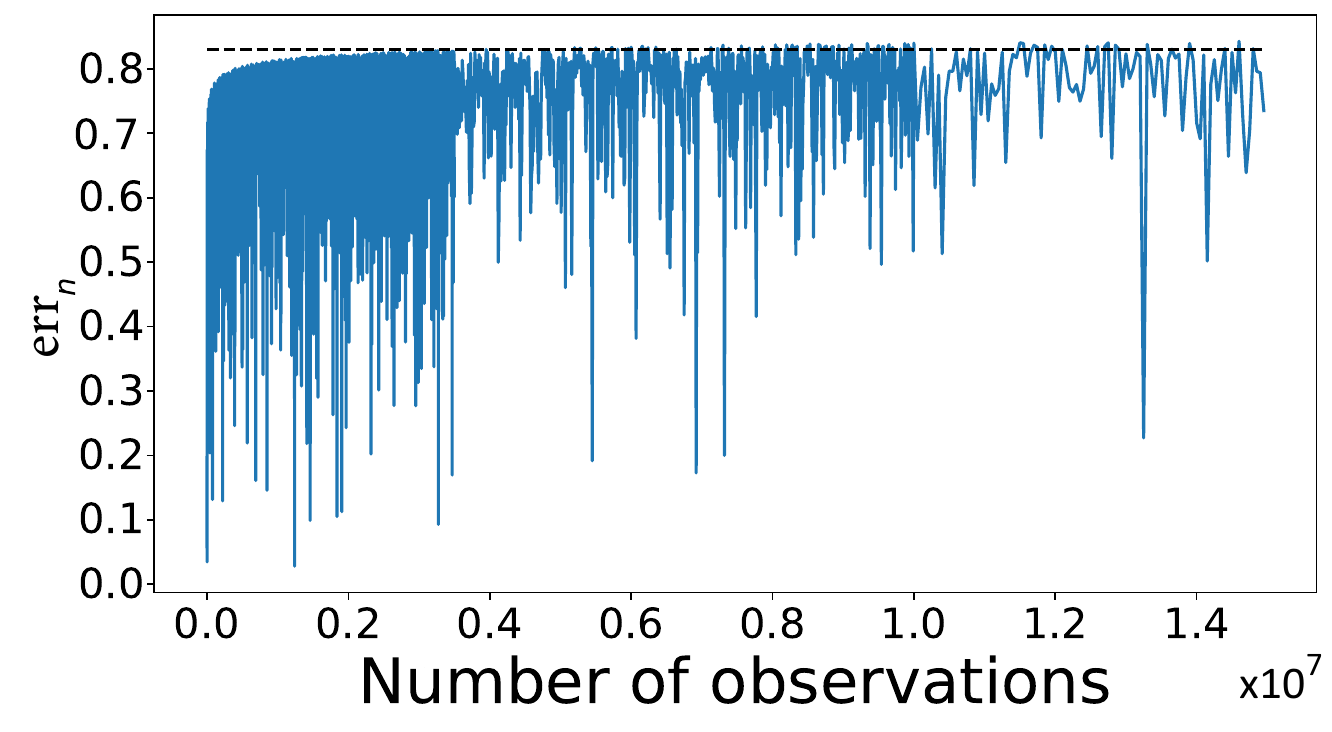}
    \caption{General structure learning for a non-maximal 
three-node graph. The true matrix $A_2^*$ is given in~(\ref{eq:matexample2}).
Left: The posterior probability $\pi(S_2^*|\Xn)$ of the true model $S_2^*$. 
Right: Fluctuations and convergence of the error term ${\rm err}_n$
in~(\ref{eq:errorn}), which, according to Theorem~\ref{thm:final2},
satisfies ${\rm err}_n\leq 1+o(1).$}
    \label{fig:general_non_maximal_case}
\end{figure*}

Next we consider observations $\{X_n\}$ generated by the 
general linear model in~(\ref{eq:LSE}), with respect to
a matrix $A_2^*$ corresponding to a {\em non-maximal} model $S_2^*$:
\begin{equation}
    \label{eq:matexample2}
    A_2^* = \begin{bmatrix}
0 & 0 & 0 \\
0 & 0 & 0 \\
1.25 & 0 & 0
\end{bmatrix}\!,
\;
S_2^* = \begin{bmatrix}
0 & 0 & 0 \\
0 & 0 & 0 \\
1 & 0 & 0
\end{bmatrix}.
\end{equation}
In view of Theorem~\ref{thm:final2}, here,
the convergence of the posterior  is expected
to be much slower, so we consider
a larger number of
$1.5\times 10^7$ observations~$\{X_n\}$.
The left plot in Figure~\ref{fig:general_non_maximal_case} 
confirms that $\pi(S_2^*|\Xn)$
converges 
much more slowly than $\pi(S_1^*|\Xn)$ did. 

In the right plot we examine the error term,
\begin{equation}
{\rm err}_n:=-\frac{2}{\log n}\log\big(1-\pi(S_2^*|\Xn)\big),
\label{eq:errorn}
\end{equation}
which, according to Theorem~\ref{thm:final2} is 
bounded above by $1+o(1)$.
Indeed, the plot   indicates that
${\rm err}_n$ fluctuates with $n$, with 
a limiting upper bound of $\approx 0.82$.

\subsection{MCMC edge detection}
\label{s:edgedet}

We compare the performance of the optimal detector 
of~\cite{NPCausalISIT}, with that of state-of-the-art methods 
for causal discovery. We use the form of the
optimal detector in Proposition~\ref{thm:OptDetector},
computed via a simple MCMC sampler as follows.
Given observations $\Xn$, in order to decide whether
$\hbchi_{ij}= 0$ or~1, we run the sampler for 
$T$ iterations, producing a collection of models
$\{S^{(t)};1\leq t\leq T\}$ 
approximately distributed according to $\pi(S|\Xn)$.
Then we average the values of 
$\bchi_{ij}(S^{(t)})$, $t=1,2,\ldots,T$ and 
declare $\hbchi_{ij}= 0$ 
if and only if that average is $\geq\tau$,
for an appropriate threshold~$\tau$.
We call this the ``optimal MCMC detector''.

Since, from~(\ref{eq:posterior_graphs_by_dimensions}),
we know $\pi(S|\Xn)$ explicitly up to normalization,
it is easy to come up with an appropriate Metropolis-Hastings
sampler. We use a simple 
irreducible proposal kernel
on $\clG_d$ which, given 
the current state $S$, proposes
a new $S'$ uniformly chosen among all neighbors of $S$,
namely, all $S'\in\clG_d$ that
differ from $S$ in exactly one edge. The remaining details are
straightforward and hence omitted. We only note that,
computationally, the most costly part
of the algorithm is the evaluation
of $\Sigma_{w}(S, \mathbf{X}^n)$, which 
requires
$O(nd^5)$ operations per MCMC iteration,
and it is the main obstacle
in scaling the algorithm  to 
high dimensions.

The performance of the optimal MCMC detector
is compared against: 
Neighborhood selection using 
LASSO~\cite{lasso}; the continuous optimization formulation 
NOTEARS~\cite{NOTEARS}; the penalized maximum likelihood (ML)
estimator in~\cite{Peters_2013}; and the two pseudo-Bayesian 
methods (``resample'' and ``bootstrap'') of~\cite{beingBayesian}. 
We consider two classes
of problems, with dimensions $d=4$ and  $d=7$.
In each case, we first randomly generate $N=1000$ $d\times d$ matrices
$A$, and for each $A$ we generate $n=10$ samples
$\Xn$ from the  model~(\ref{eq:LSE}). To generate each $A$, we 
select a model $S$ uniformly among  $S\in\clG_d$
(see, e.g.,~\cite{sprites_book_2017})
and then add i.i.d.\ $N(0,1)$ 
weights to the nonzero entries of $S$ to obtain~$A$. 

\begin{figure}[ht!]
\centering
     \includegraphics[scale=.4]{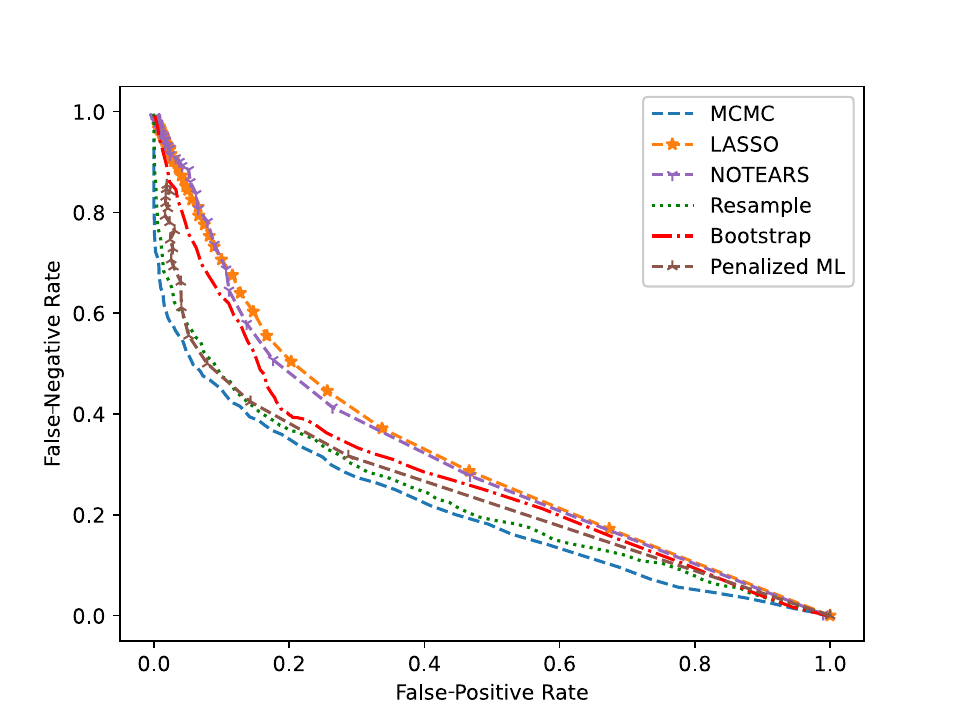}
     \includegraphics[scale=.4]{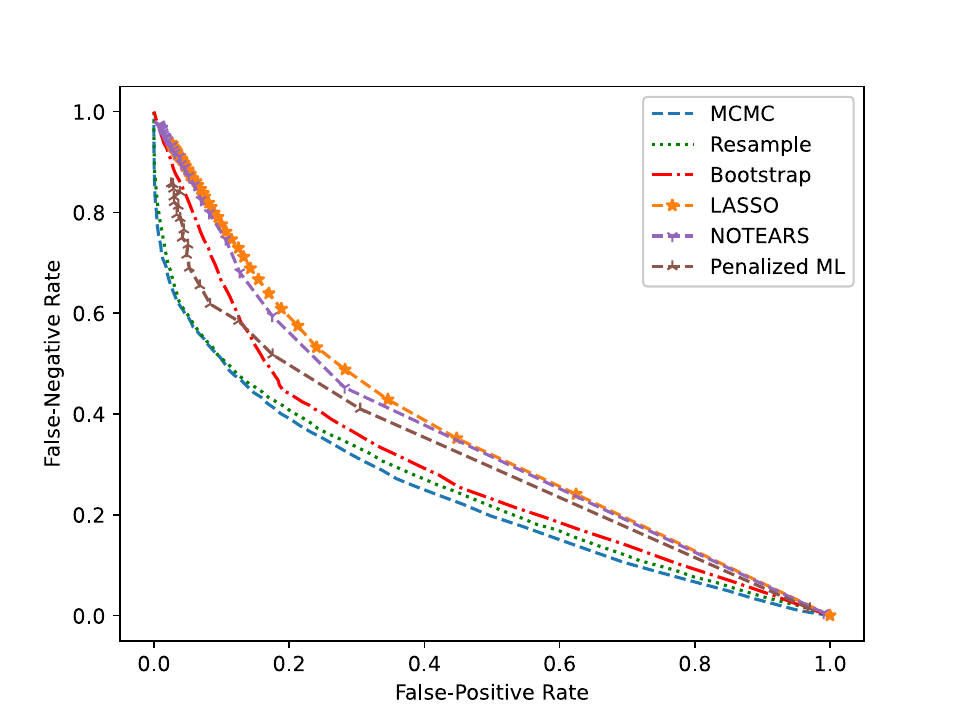}
    \caption{Performance comparison between the optimal MCMC detector
(with $T=100,000$ MCMC iterations) and state-of-the art
causal discovery methods. Top: Estimated false positive and 
false negative rates 
on simulated data from 
the linear causal model with respect to randomly chosen matrices
$A$ on $d=4$ variables. 
Bottom: Corresponding estimates 
on  simulated data
on $d=7$ variables.}
    \label{fig:samplercomparisons}
\end{figure}

For each such data set $\Xn$, each of these six methods produces 
an estimate $\hat{\bchi}$ of the support $\bchi$ of the true 
underlying model. Using these, 
we compute empirical approximations
to the false positive 
and false negative rates,
$\varepsilon^+_n$ and $\varepsilon^-_n$,
in~\eqref{def::falsepositive} and~\eqref{def::falsenegative}, respectively,
by replacing the expectations with empirical averages.
The results are shown in Figure~\ref{fig:samplercomparisons}.

The plots in Figure~\ref{fig:samplercomparisons} clearly indicate
that the optimal MCMC detector
outperforms all other methods.
The second most accurate method in almost all cases 
is the ``resample'' algorithm of~\cite{beingBayesian}.
Interestingly, the least accurate methods 
(LASSO, NOTEARS, and  penalized ML) 
 provide a single estimate $\hat{S}$ of the 
true underlying graph structure, while
the optimal MCMC detector
and the two pseudo-Bayesian methods in~\cite{beingBayesian},
actually provide probability distributions on 
$\clG_d$. When these  are close to the
actual posterior distribution
on $\clG_d$, Proposition~\ref{thm:OptDetector}
tells us that the resulting detector is near-optimal.
This explains, perhaps, the superior overall performance
of the optimal MCMC detector, as well as the fact that
the ``bootstrap''  
method of~\cite{beingBayesian} generally outperforms both LASSO 
and NOTEARS.


\appendix
\section{Appendix}

\subsection{Proofs of Lemmas~\ref{lem:relent} and~\ref{lem:minD}}
\label{appendix:BSC}

\noindent
{\sc Proof of Lemma \ref{lem:relent}.} 
Recall that,
\begin{align*}
    &D\big(N(0,\Sigma_0)\|N(0,\Sigma_1)\big)
\\&=\frac{1}{2}\left[{\rm tr}(\Sigma_1^{-1}\Sigma_0)+\log\frac{|\Sigma_1|}{|\Sigma_0|}-d\right].
\end{align*}
Here, we have $\Sigma_0=\sigma^2(I-S')^{-1}(I-S')^{-\transpose}$ and
$\Sigma_1=\sigma^2(I-S)^{-1}(I-S)^{-\transpose}$. Since 
$S^\transpose$ and $S'^\transpose$ are DAG adjacency matrices,
they are nilpotent, so the only eigenvalue of $(I-S')$ 
and $(I-S)$ is~1, therefore
the determinants $|\Sigma_0|$ and $|\Sigma_1|$ are both equal 
to $\sigma^{2d}$ and,
$$D\big(N(0,\Sigma_0)\|N(0,\Sigma_1)\big)
=\frac{1}{2}\left[{\rm tr}(\Sigma_1^{-1}\Sigma_0)-d\right].$$
Also,
$$\Sigma_1^{-1}\Sigma_0 = (I-S)^\transpose(I-S)(I-S')^{-1}(I-S')^{-\transpose},$$
and since the trace of the product of two matrices is the sum of the elements of their Hadamard product,
the trace of $\Sigma_1^{-1}\Sigma_0)$ is equal to the trace of
the product of $(I-S)(I-S')^{-1}$ and its transpose.
By one of the equivalent expressions for the Frobenius norm,
this is exactly $\|(I-S)(I-S')^{-1}\|_F^2$,
and the claimed result follows.
\qed

\smallskip

\noindent
{\sc Proof of Lemma \ref{lem:minD}.} In view of Lemma~\ref{lem:relent},
the second statement follows trivially from the first. 

For any $G\in\clG_d$, 
$(I-G)^{-1}=\sum_{k=0}^{d-1}G^k$,
since $G^k=0$ for $k\geq d$ because $G$ 
is a DAG.
Therefore, $(I-S)(I-G)^{-1}$ only has integer entries, and hence the square 
of its Frobenius norm is a nonnegative integer. Moreover, by the nonnegativity
of relative entropy and identifiability,
it has to be at least $d$. And since $P_S\neq P_G$ by 
identifiability,
it must be strictly greater than $d$, i.e., at least $d+1$.
Therefore, in order to establish the lemma, it suffices to show that, for any $G\in\clG_d$,
there is an $S\neq G$ such that:
\begin{equation}
\|(I-S)(I-G)^{-1}\|^2_F=d+1.
\label{eq:target}
\end{equation}

If $G$ is the all-zero matrix, then taking $S$ to be any 
matrix with only one~1
immediately yields~(\ref{eq:target}).
Suppose $G$ is not all-zero.
Then we can always find a pair of nodes $s\neq t$ such 
that has $s\to t$ is an edge
in $G$ (i.e., $G_{ts}=1$) and $s$ has
no parents in $G$.  Let $S\in\clG_d$ be the same as $G$ 
but with the edge $s\to t$ removed
(i.e., with $S_{ts}=0$),
and let $\Delta\in\clG_d$ consist of only that one edge, 
so that $G=S+\Delta$. Then,
\begin{align*}
    &(I-S)(I-G)^{-1}=(I-G+\Delta)(I-G)^{-1}\\&=I+\Delta(I-G)^{-1}=I+\sum_{k=0}^{d-1}\Delta G^k.
\end{align*}
Since the only nonzero entry of $\Delta$ is $\Delta_{ts}=1$, 
the product $\Delta G^k$ is all zero except 
for the $t$th row, which is 
$((G^k)_{s1},\ldots,(G^k)_{sd})$. But since $s$ has 
no parents in $G$, this entire
row must be zero for $k\geq 1$. Therefore, the only nonzero term in the sum
above is $k=0$, thus,
$$(I-S)(I-G)^{-1}=I+\Delta,$$
which implies~(\ref{eq:target})
and completes the proof.
\qed

\subsection{Proofs of Propositions~\ref{prop:exponent_approx_A}
and~\ref{prop:exponent_approx}}
\label{app:GLS}

We begin by recalling the 
classical least-squares projection formula;
see, e.g.,~\cite{hastie:book}.

\begin{lemma}
\label{lemma:regression_res}
Let $X$ be a $d$-dimensional random
vector with zero mean and positive definite
covariance $\BBE[XX^\transpose]$.
If $Y$ is an arbitrary random variable
with finite variance, then,
\begin{align*}
    &\min_{\beta \in \mathbb{R}^d} \mathbb{E} \big( Y - \beta^{\transpose}X \big)^2 \\
    &= \mathbb{E}(Y^2) - \mathbb{E}[YX^{\transpose}]\mathbb{E}[XX^T]^{-1}\mathbb{E}[YX],
\end{align*}
    and the minimum is achieved by:
    $$\beta^* = \mathbb{E}[XX^T]^{-1}\mathbb{E}[YX].$$ 
\end{lemma}

Let $\|B\|_{\rm op}$ denote the $L^2$ operator norm,
or spectral norm, of a matrix $B$.
The following is a simple bound on the operator
norm of matrix inverses.

\begin{lemma}
\label{lemma:Neumann_series_approx}
Let $Q$ be an arbitrary square matrix with
$\|Q\|_{\rm op}<1$. If $I-Q$ is invertible, then
\begin{align}
 \|(I-Q)^{-1}\|_{\rm op} &\leq \frac{1}{1-\|Q\|_{\rm op}},\nonumber\\
\mbox{and}\quad
\| (I-Q)^{-1} - 
I  \|_{\rm op} &\leq \frac{\|Q\|_{\rm op}}{1-\| Q\|_{\rm op}}.
\label{eq:opnorm}
\end{align}
\end{lemma}

\noindent 
{\sc Proof.} Since $\|Q\|_{\rm op} < 1$, the series for $(I-Q)^{-1}$ 
converges,
$(I-Q)^{-1} = \sum_{i=1}^\infty Q^i.$
For fixed $m\geq 1$, 
$$\Bigg{\|}\sum_{i=0}^m Q^i\Bigg{\|}_{\rm op} \leq \sum_{i=0}^m\|Q^i\|_{\rm op}
\leq \sum_{i=0}^m\|Q\|^i_{\rm op} \leq \sum_{i=0}^\infty\|Q\|^i_{\rm op}$$
and letting $m\to\infty$ gives the first bound.
For the second bound, we write,
$$(I-Q)^{-1} - I = \sum_{i=1}^\infty Q^i 
= Q \sum_{i=0}^{\infty}Q^i = Q(I-Q)^{-1},$$ 
so that,
$$\|(I-Q)^{-1} - I \|_{\rm op}  =\|Q(I-Q)^{-1}\|_{\rm op} 
\leq \frac{\|Q\|_{\rm op}}{1-\|Q\|_{\rm op}},$$
by the submultiplicativity of the operator norm.
\qed 

\smallskip

\noindent
{\sc Proof of Proposition~\ref{prop:exponent_approx_A}.} 
Let $1\leq j\leq d$ arbitrary. 
We first consider $\Sigma_w^{(j)}(S,\Xn)$.
From its definition~(\ref{eq:Sigmawdef}),
$\Sigma_{w}^{(j)}(S, \Xn)$ equals:
$$\Big(\frac{1}{\sigma_w^2}I+ \frac{n}{\sigma^2} 
\mathbb{E}[P_j(S, X)P_j(S, X)^{\transpose}] + L^{(j)}_n(S)\Big)^{-1},
$$
where $L_n^{(j)}(S)$ is the partial sum of centered
i.i.d.\ random matrices,
\begin{align*}
    &L^{(j)}_n(S):=\frac{1}{\sigma^2}\sum_{i=1}^n\big\{
P_j(S, X_i)P_j(S, X_i)^{\transpose}\\&\qquad \qquad\qquad\qquad \qquad-\BBE[P_j(S, X)
P_j(S, X)^{\transpose}]\big\}.
\end{align*}
We rewrite $\Sigma_w^{(j)}(S,\Xn)$ as,
\begin{equation}
\begin{aligned}
    & \frac{\sigma^2}{n}\Big[ I + \frac{\sigma^2}{n}\mathbb{E}
[P_j(S, X)P_j(S, X)^{\transpose}]^{-1}
\\&\;\;
\Big(\frac{1}{\sigma_w^2}I+L^{(j)}_n(S)\Big) \Big]^{-1} 
\mathbb{E}[P_j(S, X)P_j(S, X)^{\transpose}]^{-1}.
\end{aligned}
\label{eq:Sigmaexpand}
\end{equation}
Considering $L_n^{(j)}(S)$ as a random vector,
the multidimensional form of the LIL~\cite{ledoux} gives,
$$\|L_n^{(j)}(S)\|_2 = O(\sqrt{n \log \log n}) \quad \mbox{a.s., as}\;
n\to \infty. $$
Since the $L^2$ norm of $L_n^{(j)}(S)$ viewed as a vector is
the same as the Frobenius norm of $L_n^{(j)}(S)$ viewed as a
matrix, and since the operator norm is always bounded by
the Frobenius norm, we have that
$\|L^{(j)}_n(S)\|_{\rm op} = O(\sqrt{n \log\log n})$ a.s.,
and hence, writing $q_n(j,S)$ for
$$
\Big{\|} \frac{\sigma^2}{n}\mathbb{E}[P_j(S, X)P_j(S, X)^{\transpose}]^{-1}
\Big(\frac{1}{\sigma_w^2}I+L^{(j)}_n(S)\Big) \Big{\|}_{\rm op},
$$
we have that $q_n(j,S)=
O(n^{-1/2}\sqrt{\log\log n})$ a.s.

In particular, $q_n(j,S)<1$ for sufficiently large $n$ a.s.,
therefore,
applying~(\ref{eq:opnorm}) in 
Lemma~\ref{lemma:Neumann_series_approx} to~(\ref{eq:Sigmaexpand}),
we obtain that, eventually a.s., $\Sigma_{w}^{(j)}(S, \Xn)$ equals, 
\begin{equation}\frac{\sigma^2}{n}\mathbb{E}[P_j(S, X)P_j(S, X)^{\transpose}]^{-1} 
+ R_n(j,S),
    \label{eq:Sigma_w_approximation}
\end{equation}
where the matrix $R_n(j,S)$ satisfies,
\begin{equation}
\begin{aligned}
&\|R_n(j,S)\|_{\rm op} 
\\&\leq \Big(\frac{\sigma^2}{n}\Big)\frac{q_n(j,S)}{1-q_n(j,S)}
= O\Big( \frac{\sqrt{\log \log n}}{n^{3/2}} \Big)  
\quad\mbox{a.s.}
\end{aligned}
\label{eq:Rn}
\end{equation}
Combining this with~(\ref{eq:Sigma_w_approximation})
proves the first part of the proposition.

The second part, on the asymptotic behaviour
of $T^{(j)}(S,\Xn)$, is similar. We begin by writing,
\begin{equation}
   \begin{aligned}
 b^{(j)}(S, \Xn) 
= \frac{1}{\sigma^2} \sum_{i=1}^n X_i(j)P_j(S, X_i) 
\\= \frac{1}{\sigma^2}\Big( n\mathbb{E}[X{(j)}P_j(S, X)] 
+ \ell_n^{(j)}(S)  \Big),
   \end{aligned}
    \label{eq:b_approximation}
\end{equation}
where $\ell_n(S)=(\ell_n^{(1)}(S),\ldots,\ell_n^{(d)}(S))$ is 
the partial sum of centered
i.i.d.\ random vectors, where,
$$\ell^{(j)}_n(S) := \sum_{i = 1}^n \big\{X_i^{(j)}P_j(S, X_i) 
- \mathbb{E}[X(j)P_j(S, X)] \big\}.$$

Using the
vector form of the LIL once again, $\ell^{(j)}_n(S)$
satisfies
$\|\ell^{(j)}_n(S)\|_2 = O(\sqrt{n\log\log n})$
a.s.
Combining~(\ref{eq:Sigma_w_approximation}) and~(\ref{eq:b_approximation})
with the definition of $T^{(j)}(S,\Xn)$ in~(\ref{eq:Tj})
we see that
\begin{align*}
&T^{(j)}(S,\Xn)= \frac{1}{\sigma^4}
\Big( n\mathbb{E}[X{(j)}P_j(S, X)]^\transpose 
	+ \ell_n^{(j)}(S)^\transpose  \Big) \\
	&\qquad \Big(  \frac{\sigma^2}{n}
	\mathbb{E}[P_j(S, X)P_j(S, X)^{\transpose}]^{-1} + R_n(j,S)\Big) 
	\\
&\qquad \qquad \qquad
	\Big( n\mathbb{E}[X{(j)}P_j(S, X)] + \ell^{(j)}_n(S)  \Big),
\end{align*}
which, by~(\ref{eq:Rn}), gives, 
$$
T^{(j)}(S,\Xn)
=n T_\infty^{(j)}(S) + O(\sqrt{n\log\log n}) \quad \mbox{a.s.}$$
Finally, the alternative form of $T_\infty^{(j)}(S)$ in the proposition
follows by Lemma~\ref{lemma:regression_res}.
\qed

\smallskip

\noindent
{\sc Proof of Proposition~\ref{prop:exponent_approx}.} 
Recalling the definition of $\mu_w^{(j)}(S,\Xn)$ in~(\ref{eq:muwj}) and
using~(\ref{eq:Sigma_w_approximation})
and~(\ref{eq:b_approximation}), we have,
\begin{align*}
    &\mu_w^{(j)}(S, \Xn) \\
&= \Big(\frac{\sigma^2}{n}\mathbb{E}[P_j(S, X)P_j(S, X)^{\transpose}]^{-1} 
+ R_n(j,S)\Big) \\& \qquad\qquad \frac{1}{\sigma^2}
\Big( n\mathbb{E}[X{(j)}P_j(S, X)] + \ell_n^{(j)}(S)  \Big)\\
&= \mu_\infty^{(j)}(S) + r_n\quad\mbox{a.s.},
\end{align*}
with $\|r_n\|_2 \leq O(\sqrt{n^{-1}\log\log n})$ a.s.
Again, the alternative form of $\mu_\infty^{(j)}(S)$ 
follows by Lemma~\ref{lemma:regression_res}.

For $S^*$, $\mu^{(j)}_\infty(S^*)$
equals,
\begin{align*}
    &\mathbb{E}[P_j(S^*, X)P_j(S^*, X)^\transpose]^{-1} \mathbb{E}[P_j(S^*, X)X{(j)}] \\
    & = \mathbb{E}[P_j(S^*, X)P_j(S^*, X)^\transpose]^{-1} \\& \qquad \qquad \quad  \mathbb{E}[P_j(S^*, X)(P_j(S^*, X)^\transpose\mathbf{w}^*{(j)} + \epsilon{(j)})] \\
    & = \mathbf{w}^*{(j)} \\ &\quad  + \mathbb{E}[P_j(S^*, X)P_j(S^*, X)^\transpose]^{-1} \mathbb{E}[P_j(S^*, X)\epsilon{(j)}]\\
    & = \mathbf{w}^*{(j)},
\end{align*} 
where the last equality 
follows from the definition of $P_i(S^*, X)$ 
and the fact that $S^*$ is DAG,
which imply that $\mathbb{E}[P_j(S^*, X)\epsilon{(j)}] = 0$
for all $j$.
\qed

\subsection{Proofs of Lemmas~\ref{lemma:maximal_case}
and~\ref{lemma:non_maximal_case}} 
\label{app:mainlemmas}

\noindent
{\sc Proof of Lemma~\ref{lemma:maximal_case}.} 
From~(\ref{eq:posterior_graphs_by_dimensions}), we have,
\begin{align*}
    &\log\frac{\pi(S|\mathbf{X}^n)}{\pi(S^*|\mathbf{X}^n)}  
\\&= \sum_{j=1}^d \Big\{\frac{1}{2}\big(T^{(j)}(S, \Xn)-T^{(j)}(S^*, \Xn)\big) 
\\& \quad+ \frac{1}{2}\big(\log|\Sigma_{w}^{(j)}(S,\Xn)| - \log|\Sigma_{w}^{(j)}(S^*,\Xn)| \big) \Big\}.
\end{align*}
We first examine the first term in the exponent.
Using Propositions~\ref{prop:exponent_approx_A} 
and~\ref{prop:exponent_approx},
as $n\to\infty$, we have,
\begin{align*}
&T^{(j)}(S, \Xn) - T^{(j)}(S^*, \Xn)
\\ &= nT_\infty^{(j)}(S) - n T_\infty^{(j)}(S^*)+O(\sqrt{n\log\log n})  \\
    &  
\overset{(a)}{=} \frac{n}{\sigma^2} \Big(
\mathbb{E}\Big[\big(X{(j)}-\mathbf{w}{(j)}^{*\transpose}
	P_j(S^*, X)\big)^2\Big]\\
&\quad
 -\mathbb{E}\Big[\big(X{(j)}
	-\mu_\infty^{(j)}(S)^\transpose P_j(S, X)\big)^2\Big]\Big) \\
&\quad
	+O(\sqrt{n\log\log n})\\
    &=  \frac{n}{\sigma^2}\mathbb{E}[\epsilon{(j)^2}] \!-\! \frac{n}{\sigma^2}
\mathbb{E}\Big[\big(X{(j)}-\mu_\infty^{(j)}(S)^\transpose P_j(S, X)\big)^2\Big] \\&\quad 
+ O(\sqrt{n\log\log n})\\
    &= n - \frac{n}{\sigma^2}
\mathbb{E}\Big[\big(X{(j)}-\mu_\infty^{(j)}(S)^\transpose P_j(S, X)\big)^2
\Big]\\ & \quad +O(\sqrt{n\log\log n}),\quad\mbox{a.s.,}
\end{align*}
where $(a)$ follows from~(\ref{eq:def_T_i}) and 
the fact that,
by Proposition~\ref{prop:exponent_approx},
$\mu^{(j)}_\infty(S^*) = \mathbf{w}{(j)}^*$. 

Similarly, for the second term, in the notation
of the proof of Proposition~\ref{prop:exponent_approx_A},
for any $S\neq S^*$, not necessarily assuming
$S^*$ is a subgraph of $S$, we have 
\begin{align*}
&\frac{1}{2} \log \big| \Sigma_{w}^{(j)}(S, \Xn)\big| 
- \frac{1}{2} \log \big| \Sigma_{w}^{(j)}(S^*, \Xn)\big|\\ 
&= \frac{1}{2} \log \Big{|} \frac{1}{n}\mathbb{E}[P_j(S, X)P_j(S, X)^{\transpose}]^{-1}\! + R_n(j,S)\Big{|} \nonumber \\ &- \frac{1}{2} \log \Big{|} \frac{1}{n}\mathbb{E}[P_j(S^*, X)P_j(S^*, X)^{\transpose}]^{-1} \!+ R_n(j,S^*)\Big{|}.
\end{align*}
Hence, a.s.\ as $n\to\infty$,
\begin{align}
&\frac{1}{2} \log \big| \Sigma_{w}^{(j)}(S, \Xn)\big| 
- \frac{1}{2} \log \big| \Sigma_{w}^{(j)}(S^*, \Xn)\big| 
\nonumber\\
& = \frac{1}{2}\big[\|P_j(S, X) \|_0 - \|P_j(S^*, X) \|_0\big]\log n \nonumber
\\ &\quad+O\big(n^{-3/2}\sqrt{\log \log n}\big),
\label{eq:Sigmaexp}
\end{align}
where $\|a\|_0$ is the number of nonzero components
of~$a$. 

Substituting these two estimates to the expression
for the log-odds, we have, almost surely as $n\to\infty$,
\begin{align*}
	&\log\frac{\pi(S|\mathbf{X}^n)}{\pi(S^*|\mathbf{X}^n)}
\\ &= 
	\frac{n}{2\sigma^2}\sum_{i=1}^d \!
	\Big\{\sigma^2 \!-  \mathbb{E}\Big[\big(X{(j)}
	-\mu_\infty^{(j)}(S)^\transpose P_j(S, X)\big)^2\Big]
	\!\Big\}\\
&\qquad
	 + \Big(\frac{\log n}{2}\Big) 
	\sum_{i =1}^d [\|P_i(S, X) \|_0 
	- \|P_i(S^*, X) \|_0] \\ &\qquad \qquad\qquad \qquad \qquad\qquad \quad+O(\sqrt{n\log \log n}),
\end{align*}
and, letting $A := m(S, \mu_\infty(S))$,
\begin{align*}
	&\log\frac{\pi(S|\mathbf{X}^n)}{\pi(S^*|\mathbf{X}^n)}
= 
	\frac{n}{2\sigma^2}\Big(d\sigma^2 - \mathbb{E}\big[\|X-AX\|_2^2
	\big]\Big) \\& \qquad+ 
	\frac{\log n}{2} [\|S^*\|_0-\|S\|_0]+O(\sqrt{n\log \log n})\\
& 
	= \frac{n}{2}\Big(d 
	- \mathbb{E}[\|(I-A)X\|_2^2]\Big) +O(\sqrt{n\log \log n})\\
  \notag  & = \frac{n}{2}\Big(d- \|(I-A)(I-A^*)^{-1}\|_F^2\Big) \\ &\qquad \qquad \qquad \qquad \qquad \qquad+O(\sqrt{n\log \log n})\\
   & = -nD(P_{A^*} \| P_{A}) +O(\sqrt{n\log \log n}) \quad \mbox{a.s.},
\end{align*}
where the last step follows from Lemma~\ref{lem:relent2}.
Finally, 
since $S^*$ is not a subgraph of $S$, we have
$A \neq A^*$, $P_{A^*}\neq P_A$ and 
$D(P_{A^*}\|P_A) > 0$
by identifiability.
\qed

\smallskip

\noindent
{\sc Proof of Lemma~\ref{lemma:non_maximal_case}.} 
Since $S^+$ and $S^*$ differ in only one edge,
we may assume
without loss of generality that
$P_j(S^+,X)=P_j(S^*,X)$ for $j\geq 2$,
while $P_1(S^+,X) =[P_1(S^*, X)^\transpose, X(k)]^\transpose$
for a fixed $k\in\{1,\ldots,d\}$ different for the indices
$j$ of the $X(j)$'s that are present in $P_1(S^*,X)$.

Our starting point is the expression for the posterior
log-odds in the beginning of the proof of Lemma~\ref{lemma:maximal_case}.
By our assumptions,
for all $j\neq 1$ we have
$T^{(j)}(S^+,\Xn)=T^{(j)}(S^*,\Xn)$
and $\Sigma_w^{(j)}(S^+,\Xn)=\Sigma_w^{(j)}(S^*,\Xn)$.
Therefore,
\begin{align}
&    \log\frac{\pi(S^+|\mathbf{X}^n)}{\pi(S^*|\mathbf{X}^n)}  
= \frac{1}{2}\big(T^{(1)}(S^+, \Xn)-T^{(1)}(S^*, \Xn)\big) 
\nonumber\\
&+\! \frac{1}{2}\big(\log|\Sigma_{w}^{(1)}(S^+,\Xn)| 
\!-\! \log|\Sigma_{w}^{(1)}(S^*,\Xn)| \big).
\label{eq:starting}
\end{align}
Most of the proof will be devoted
to the evaluation of the first term in
the right-hand side of~(\ref{eq:starting}).

Recalling the definition of 
$T^{(j)}(S^+, \Xn)$ in~(\ref{eq:Tj}),
we see that 
$T^{(1)}(S^+, \Xn)$ can be written as,
    \begin{align}
        &T^{(1)}(S^+, \Xn) 
    = \begin{bmatrix}
b^{(1)}(S^*, \Xn)\\
\frac{1}{\sigma^2}\sum_{i = 1}^n X_i(k)X_i(1)
\end{bmatrix}^\transpose \nonumber\\& \;
\bar{\Sigma}(S^+,S^*,\Xn)^{-1}
\begin{bmatrix}
b^{(1)}(S^*, \Xn)\\
\frac{1}{\sigma^2}\sum_{i = 1}^n X_i(k)X_i(1)
\end{bmatrix},
\label{eq:Tn1}
\end{align}
where, $\bar{\Sigma}(S^+,S^*,\Xn)$ is given by \eqref{eq:Tn2} at the bottom of the page.
\begin{figure*}[b]
\hrulefill 
    \begin{equation}
\bar{\Sigma}(S^+,S^*,\Xn)
:=\begin{bmatrix}
\frac{1}{\sigma_w^2}I + \frac{1}{\sigma^2}\sum\limits_{i=1}^n P_1(S^*, X_i) P_1(S^*, X_i)^\transpose & \frac{1}{\sigma^2}\sum\limits_{i = 1}^n 
X_i(k)P_1(S^*, X_i) \\[1em]
\frac{1}{\sigma^2}\sum\limits_{i = 1}^n X_i(k) P_1(S^*, X_i)^\transpose & \frac{1}{\sigma_w^2} + \frac{1}{\sigma^2}\sum\limits_{i = 1}^n X_i(k)^2
\end{bmatrix}
\label{eq:Tn2}
\end{equation}
\begin{equation}
\bar{\Sigma}(S^+,S^*,\Xn)^{-1} \!
= \begin{bmatrix}
\Sigma_w^{(1)}(S^+, \Xn) + a_n^{-1}\Sigma_w^{(1)}(S^+, \Xn) d_n d_n^T \Sigma_w^{(1)}(S^+, \Xn) & -a_n^{-1}\Sigma_w^{(1)}(S^+, \Xn) d_n \\
-a_n^{-1} d_n^T \Sigma_w^{(1)}(S^+, \Xn) & a_n^{-1}
\end{bmatrix}
\label{eq:Tn3}
\end{equation}
\end{figure*}
By the 
Sherman-Morrison inversion formula~\cite{sherman},
$\bar{\Sigma}(S^+,S^*,\Xn)^{-1}$ is given by \eqref{eq:Tn3} at the bottom of the page
where,
\begin{align*}
a_n & 
:= 1+\sum_{i = 1}^n X_i(k)^2 - d_n^\transpose\Sigma_w^{(1)}(S^+,\Xn)d_n,\\
d_n &
:= \frac{1}{\sigma^2}\sum_{i = 1}^n X_i(k) P_1(S^*,X_i).
\end{align*}
Since $S^*$ is the true underlying structure,
we have,
$$X_i{(1)} = P_1(S^*, X_i)^\transpose\mathbf{w}^*{(1)}  + \epsilon_i{(1)},$$
where, from~(\ref{eq:def_T_i}), $\ww^*(1)$ equals
$$ \mathbb{E}[P_1(S^*, X)P_1(S^*, X)^\transpose]^{-1}
\mathbb{E}[X{(i)}P_1(S^*, X)].$$ 
Therefore,
\begin{align*}
&\begin{bmatrix}
b^{(1)}(S^*, \Xn)\\
\frac{1}{\sigma^2}\sum_{i = 1}^n X_i(1)X_i(k)
\end{bmatrix} 
\\&= 
\frac{1}{\sigma^2}\!\begin{bmatrix}
\sum_{i=1}^n P_1(S^*, X_i) (P_1(S^*, X)^\transpose\mathbf{w}^*{(1)}\! +\! \epsilon_i{(1)})\\
\sum_{i = 1}^n X_i(k)[\mathbf{w}^*(1) P_1(S^*, X_i) + \epsilon_i{(1)}]
\end{bmatrix}\\
&= 
\frac{1}{\sigma^2}\begin{bmatrix}
\big[\big(\Sigma_w^{(1)}\big)^{-1}-\frac{1}{\sigma_w^2}I\big] 
\mathbf{w}^*{(1)} + c_n\\
\mathbf{w}^*{(1)} d_n + e_n
\end{bmatrix},
\end{align*}
where,
\begin{align*}
c_n &:= \sum_{i=1}^n \epsilon_i{(1)} P_1(S^*, X_i), \quad
e_n := \sum_{i=1}^n\epsilon_i{(1)}X_i(k).
\end{align*}

Now we turn to asymptotics.
For $d_n$, using the multidimensional
LIL again, we have
that, as $n\to\infty$,
\begin{equation}
    d_n = n\mathbb{E}[X(k)P_1(S^*, X)] + v_n,
    \label{eq:approx_d_n}
\end{equation}
where $\|v_n\|_2=O(\sqrt{n\log\log n})$ almost surely.
For $\Sigma^{(1)}_w(S^+, \Xn)$ we have,
from~(\ref{eq:Sigma_w_approximation}),
a.s.\ as $n\to\infty$,
\begin{align}
\Sigma_{w}^{(1)}(S^+, \Xn) 
&= \frac{1}{n}\mathbb{E}[P_1(S^+, X)P_1(S^+, X)^{\transpose}]^{-1}\nonumber\\ 
&\quad +R_n(1,S^+),
    \label{eq:approx_inverse}
    \end{align}
where $\|R_n(1,S^+)\|_{\rm op}=
O(n^{-3/2}\sqrt{\log \log n})$.
And for $a_n$, combining~(\ref{eq:approx_d_n})
and~(\ref{eq:approx_inverse}),
we obtain that,  
\begin{align}
&a_n 
= n\mathbb{E}[X(k)^2] \nonumber \\ \nonumber&- n\mathbb{E}[X(k)P_1(S^*, X)^\transpose]
	\mathbb{E}[P_1(S^*, X)P_1(S^*, X)^\transpose]^{-1}
	\\  &\quad \mathbb{E}[X(k)P_1(S^*, X)]
	 +O(\sqrt{n \log \log n})\quad\mbox{a.s.}
    \label{eq:approx_alpha}
\end{align}
It is easy to see that, by
Lemma~\ref{lemma:regression_res}, the coefficient of the linear 
term in~(\ref{eq:approx_alpha}) is nonnegative:
\begin{align*}
&\mathbb{E}[X(k)^2] 
\\&- \mathbb{E}[X(k)P_1(S^*, X)^\transpose]\mathbb{E}
[P_1(S^*, X)P_1(S^*, X)^\transpose]^{-1}\\&\qquad \qquad\qquad\qquad\qquad\qquad\quad\mathbb{E}[X(k)P_1(S^*, X)] \\
&= \mathbb{E}[(X(k)-\beta^{*\transpose}P_1(S^*, X))^2],
\end{align*} 
with 
$$\beta^* 
= \mathbb{E}[P_1(S^*, X)P_1(S^*, X)^\transpose]^{-1}\mathbb{E}[RP_1(S^*, X)].$$
Moreover, it is  zero iff $X(k)$ is a linear 
function of $P_1(S^*, X)$, which is impossible by the causal minimality
of the model.
Therefore, $a_n>0$ eventually a.s.

Now we turn to our main task,
namely, the evaluation of the first term
in the right-hand side of~(\ref{eq:starting}).
We observe that, 
using~(\ref{eq:Tn1}),~(\ref{eq:Tn2}) and~(\ref{eq:Tn3}),
\begin{align*}
&T^{(1)} (S^+, \Xn)
\\&= 
	b^{(1)}(S^*, \Xn)^\transpose \Sigma_w^{(1)}(S^*, \Xn)b^{(1)}(S^*, \Xn)  \\
&\quad
	+a_n^{-1}b^{(1)}(S^*, \Xn)^\transpose\Sigma_w^{(1)}(S^*, \Xn)d_nd_n^\transpose \\&\qquad\qquad\qquad\qquad\qquad\Sigma_w^{(1)}(S^*, \Xn)b^{(1)}(S^*, \Xn)\\
&\quad
	-2a_n^{-1}b^{(1)}(S^*, \Xn)^\transpose\Sigma_w^{(1)}(S^*, \Xn)
	d_n\\&\qquad\qquad(\mathbf{w}^*(1)^\transpose d_n+e_n) 
	+ a_n^{-1}(\mathbf{w}^*(1)^\transpose d_n+e_n)^2,
\end{align*}
while $T^{(1)}(S^*, \Xn)$ equals:
\begin{equation*}
b^{(1)}(S^*, \Xn)^\transpose 
\Sigma_w^{(1)}(S^*, \Xn)b^{(1)}(S^*, \Xn).
\end{equation*}
Combining these two, gives,
\begin{align*}
    &T^{(1)}(S^+, \Xn) = T^{(1)}(S^*, \Xn) + a_n^{-1}
	\Big(c_n-\frac{\mathbf{w}^*{(1)}}{\sigma_w^2}\Big)^\transpose\\
    &\quad\Sigma_w^{(1)}(S^*, \Xn)d_nd_n^\transpose\Sigma_w^{(1)}(S^*, \Xn)
	\Big(c_n-\frac{\mathbf{w}^*{(1)}}{\sigma_w^2}\Big) \\
    & - 2a_n^{-1}\Big(c_n-\frac{\mathbf{w}^*{(1)}}{\sigma_w^2}\Big)^\transpose\Sigma_w^{(1)}(S^*, \Xn)d_ne_n + a_n^{-1}e_n^2,
\end{align*}
and hence,
\begin{align*}
&T^{(1)}(S^+, \Xn) - T^{(1)}(S^*, \Xn)\\& 
=a_n^{-1}+
	\Bigg( \Big(c_n - \frac{\mathbf{w}^*{(1)}}{\sigma_w^2}\Big)^\transpose\Sigma_w^{(1)}(S^*, \Xn)d_n - e_n   \Bigg)^2.
\end{align*}
Since $a_n$ is eventually a.s.\ positive,
we have:
\begin{equation}
\begin{aligned}
    \delta_n:=T^{(1)}(S^+, \Xn) - T^{(1)}(S^*, \Xn)\geq 0
\\\quad\mbox{eventually a.s.}
\end{aligned}
\label{eq:zetan}
\end{equation}

For the second term in the right-hand side of~(\ref{eq:starting}),
we observe that, as $n\to\infty$, the general expansion~(\ref{eq:Sigmaexp})
in this case reduces to,
\begin{equation}
\begin{aligned}
    &\frac{1}{2} \log \big| \Sigma_{w}^{(j)}(S^+, \Xn)\big| 
- \frac{1}{2} \log \big| \Sigma_{w}^{(j)}(S^*, \Xn)\big|\\ &= 
\frac{1}{2}\log n 
+O\big(n^{-3/2}\sqrt{\log \log n}\big)
\quad\mbox{a.s.}
\end{aligned}
\label{eq:Sigmaexp2}
\end{equation}

Substituting the bounds~(\ref{eq:zetan}) and~(\ref{eq:Sigmaexp2})
into~(\ref{eq:starting}), yields the claimed result.
\qed


\bibliographystyle{plain}
\bibliography{references}

\end{document}